\RequirePackage{fix-cm}
\documentclass[smallextended]{svjour3}       
\smartqed  
\usepackage{graphicx}
\usepackage{amsmath,amssymb,amsfonts}
\usepackage{mathrsfs}
\usepackage{mathptmx}      

\usepackage{amsmath,amssymb,amsfonts}
\usepackage{graphicx,graphics}
\usepackage{subfigure}
\usepackage{multirow}
\usepackage{pgfplots}
\usetikzlibrary{plotmarks}
\usepackage{hyperref}

\usepackage{geometry}
\usepackage{listings}
\usepackage{color}

\definecolor{dkgreen}{rgb}{0,0.6,0}
\definecolor{gray}{rgb}{0.5,0.5,0.5}
\definecolor{mauve}{rgb}{0.58,0,0.82}

\newcommand{\Eref}[1]{Equation (\ref{#1})}
\newcommand{\fref}[1]{Figure (\ref{#1})}

\newcommand{\rmd}{\mathrm{d}}
\newcommand{\bveps}{\boldsymbol{\varepsilon}}
\newcommand{\bvsig}{\boldsymbol{\sigma}}

\newcommand{\bb}{\mathbf{b}}
\newcommand{\bigb}{\mathbf{B}}

\newcommand{\dd}{\mathbf{D}}

\newcommand{\kk}{\mathbf{K}}
\newcommand{\bn}{\mathbf{N}}

\newcommand{\qq}{\mathbf{q}}

\newcommand{\uu}{\mathbf{u}}

\newcommand{\xx}{\mathbf{x}}

\begin{document}

\title{Displacement based finite element formulations over polygons: a comparison between Laplace interpolants, strain smoothing and scaled boundary polygon formulation
}

\titlerunning{Polygonal finite element formulations}        

\author{Sundararajan Natarajan         \and
        Ean Tat Ooi \and Irene Chiong \and Chongmin Song 
}


\institute{Sundararajan Natarajan \at
              School of Civil \& Environmental Engineering \\
              The University of New South Wales, Sydney, NSW 2052 \\
              Australia.
              Tel.: +61 2 9385 5030\\
              \email{sundararajan.natarajan@gmail.com}           
           \and
           Ean Tat Ooi \at
              School of Civil \& Environmental Engineering \\
              The University of New South Wales, Sydney, NSW 2052 \\
              Australia.
              \and
              Irene Chiong \at
              School of Civil \& Environmental Engineering \\
              The University of New South Wales, Sydney, NSW 2052 \\
              Australia.
              \and
               Chongmin Song\at
              School of Civil \& Environmental Engineering \\
              The University of New South Wales, Sydney, NSW 2052 \\
              Australia.
}

\date{Received: date / Accepted: date}

\maketitle

\begin{abstract}
Three different displacement based finite element formulations over arbitrary polygons are studied in this paper. The formulations considered are: the conventional polygonal finite element method (FEM) with Laplace interpolants, the cell-based smoothed polygonal FEM with simple averaging technique and the scaled boundary polygon formulation. For the purpose of numerical integration, we employ the sub-traingulation for the polygonal FEM and classical Gaussian quadrature for the smoothed FEM and for the scaled boundary polygon formulation. The accuracy and the convergence properties of these formulations are studied with a few benchmark problems in the context of linear elasticity and the linear elastic fracture mechanics. The extension of scaled boundary polygon to higher order polygons is also discussed.
\end{abstract}

\section{Introduction}
The finite element method (FEM) is a versatile technique for the numerical approximation of solutions of partial differential equations (PDEs). Traditional FEM simulations rely on strictly tetrahedral or hexadedral meshes in 3D (or triangular, quadrilateral meshes in 2D). In generating a FE mesh, a balance is required between the accuracy and the flexibility in the mesh generation. For example, triangulation of a domain is relatively easy when compared with quadragulation, whilst the quadrilateral mesh is more accurate than the triangular mesh. It is relatively easy to construct interpolants over standard shapes, viz., triangles and tetraehdrals. The use of standard shapes, viz., triangles (or quadrilaterals) and tetrahedrals (or hexahedral) simplifies the approach, however, allowing only a few element shapes can be too restrictive, because
\begin{itemize}
\item it may require sophisticated meshing algorithm to generate high-quality mehes, esp. with quadrilaterals, for meshing complex geometries;
\item it may require complex remeshing to capture topological changes, for instance due to discontinuous surface growth.
\end{itemize}

\subsection{Background}
With a goal to decrease the constraints imposed on finite element meshes different methods have been introduced. For instance, the meshfree methods~\cite{friesmatthies2003}, the partition of unity methods (PUMs)~\cite{duarteoden1996,friesbelytschko2010}, the smoothed finite element methods (SFEM)~\cite{liutrung2010,bordasnatarajan2011} and the recent, isogeometric analysis~\cite{cottrellhughes2009}. However, the PUMs and the SFEM still require the domain to be discretized using a combination of triangles and quadrilaterals in 2D, unless, the boundary is defined implicitly. This causes additional difficulties in imposing boundary conditions and ensuring an accurate definition of the boundary of the domain~\cite{belytschkoparimi2003,moumnassibelouettar2011}. The introduction of the isogeometric analysis has revolutionized the analysis procedure. It circumvents the need to discretize the domain with standard shapes and provides a natural link with the CAD models.

On another, related front, generalizations of FEM on arbitrary polygonal and polyhedral meshes have been the subject of increasing attention in the research community, both in computational physics \cite{moorthyghosh2000,dasgupta2003,szesheng2005,pavankumarjayabal2010,jayabalmenzel2011,krausrajagopal2013} and in computer graphics~\cite{floater2003,warrenschaefer2007}.

\paragraph{Polygonal FEM} In polygonal finite elements, the use of elements with more than four sides can provide flexibility, especially in meshing and accuracy~\cite{sukumartabarraei2004}. Ghosh \textit{et al.,}~\cite{moorthyghosh2000} developed the Vorono\"{i} cell finite element model (VCFEM) to model the mechanical response of heterogeneous microstructures of composites and porous materials with heterogeneities. The VCFEM is based on the assumed stress hybrid formulation and was further developed by Tiwary \textit{et al.,}~\cite{tiwaryhu2007} to study the behaviour of microstructures with irregular geometries. Rashid and Gullet~\cite{rashidgullet2000} proposed a variable element topology finite element method (VETFEM), in which the shape functions are constructed using a constrained minimization procedure. Sukumar~\cite{sukumar2003} used Vorono\"{i} cells and natural neighbor interpolants to develop a finite difference method on unstructured grids. Biabanaki and Khoei~\cite{biabanakikhoei2012} employed the polygonal FEM technique to model the large deformation response of interfaces. In their approach, the polygonal finite elements were generated from a non-conforming regular mesh. Wachspress interpolants were employed over the polygons.  Kraus \textit{et al.,}~\cite{krausrajagopal2013} presented a polygonal FEM based on the constrained adaptive Delaunay tessellation. The polygonal elements can also be used as transition elements to simplify meshing or to describe the microstructure of polycrystalline alloys~\cite{szesheng2005,pavankumarjayabal2010,jayabalmenzel2011} in a rather straightforward manner. However, approximation functions on polygonal elements are usually non-polynomial, which introduces difficulties in numerical integration. Improving numerical integration over polytopes have gained increasing attention~\cite{sukumartabarraei2004,natarajanbordas2009,mousavixiao2010,talischipaulino2013}.

\paragraph{Polygonal SFEM} In the stabilized conforming nodal integration~\cite{chenwu2001}, the strain is written as the divergence of a spatial average of the standard (compatible) strain field - i.e., the symmetric gradient of the displacement field. This concept was incorporated into the FEM by Liu~\cite{liunguyen2007} and extensively studied in~\cite{nguyen-xuanbordas2008,bordasnatarajan2011,nguyen-thoiliu2011}. Depending on the number and geometry of the subcells used, a spectrum of methods, each exhibiting a set of unique properties. Interested readers are referred to the literature~\cite{liutrung2010} and references therein. Generalizations of the SFEM to arbitrarily shaped polygons were reported in~\cite{dailiu2007,nguyen-thoiliu2011}.

\paragraph{Scaled boundary polygon}
Wolf and Song~\cite{wolfsong2001} introduced the scaled boundary finite element method (SBFEM) for elasto-statics and elasto-dynamic problems. The SBFEM reduces the governing PDE to a set of ordinary differential equations. Like the FEM, no fundamental solution is required and like the boundary element method, the spatial dimension is reduced by one, since only the boundary need to be discretized, resulting in a decrease in the total degrees of freedom. The SBFEM relies on defining a `scaling center' from which the entire domain is visible. By exploiting the unique feature of the scaling center, the method allows the computation of stress intensity factor directly from their definitions. Natarajan and Song~\cite{natarajansong2013} combined the extended FEM and the SBFEM, thus, circumventing the need to know a priori the enrichment functions, required by the former. Recently, Ooi \textit{et al.,}~\cite{ooisong2012} employed scaled boundary formulation in polygonal elements to study crack propagation, but its performance and application in the context of linear elasticity was not studied. 

\subsection{Objective}
The generalization of finite elements over polygonal and polyhedral elements is a subject of increasing attention in the research community. Apart from the aforementioned formulations, recent studies, among others include developing polygonal elements based on the virtual nodes~\cite{tangwu2009} and the virtual element methods~\cite{veigabrezzi2013}. To the author's knowledge, a comparison of different displacement based formulations over arbitrary polygons has not been reported yet in the literature. The main objective of the paper is to study the accuracy and the convergence properties of different displacement based finite element formulations, esp., the polygonal FEM, the polygonal SFEM and the scaled boundary polygon in the context of linear elasticity and linear elastic fracture mechanics. The shape functions used in these formulations are: Laplace interpolants for the polygonal FEM, simple averaging technique for the polygonal nSFEM and Gauss-Legendre shape functions for the scaled boundary polygon formulation. The scaled boundary polygons with higher order shape functions is also studied.

\subsection{Outline} The paper is organized as follows. Section \ref{fempoly} presents the governing equations and weak form for 2D static elasticity problem. The different displacement based finite element formulations considered in this study are also discussed. In Section \ref{numexamp}, with a few benchmark problems from linear elasticity, the accuracy and the convergence properties of various finite element techniques are studied. Later, the scaled boundary polygon formulation is applied to study problems in linear elastic fracture mechanics. The results from scaled boundary polygon formulation are compared with results available in the literature, for example with the extended FEM, followed by concluding remarks in the last section.


\section{Overview of finite element techniques over polygons}
\label{fempoly}
\subsection{Governing equations and weak form}
For a 2D static elasticity problem defined in the domain $\Omega$ bounded by $\Gamma = \Gamma_u \bigcup \Gamma_t$, $\Gamma_u \bigcap \Gamma_t = \emptyset$, in the absence of body forces, the governing equation is given by:
\begin{equation}
\nabla_s \cdot \bvsig = \mathbf{0} \hspace{0.25cm} \textup{in} \hspace{0.25cm}  \Omega
\end{equation}
with the following conditions prescribed on the boundary: 
\begin{align}
\uu &= \overline{\uu} \hspace{0.25cm} \textup{in} \hspace{0.25cm}  \Gamma_u \nonumber \\
\bvsig \cdot \mathbf{n} &= \overline{\mathbf{t}} \hspace{0.25cm} \textup{on} \hspace{0.25cm}  \Gamma_t
\end{align}
where $\bvsig$ is the stress tensor. The discrete equations for this problem are formulated from the Galerkin weak form:
\begin{equation}
\int\limits_{\Omega} \left( \nabla_s \uu \right)^{\rm T} \dd \left( \nabla_s \delta \uu \right)~\rmd \Omega - \int\limits_\Omega \left( \delta \uu \right)^{\rm T} \bb~\rmd \Omega - \int\limits_{\Gamma_t} \left( \delta \uu \right)^{\rm T} \overline{\mathbf{t}}~\rmd \Gamma = \mathbf{0}
\label{eqn:weakform}
\end{equation}
where $\uu$ and $\delta \uu$ are the trial and the test functions, respectively and $\dd$ is the material constitutive matrix. The FEM uses the following trial and test functions:
\begin{equation}
\uu^h(\xx) = \sum\limits_{I=1}^{NP} \bn_I(\xx) \mathbf{d}_I, \hspace{0.25cm} \delta\uu^h(\xx) = \sum\limits_{I=1}^{NP} \bn_I(\xx) \delta \mathbf{d}_I
\label{eqn:trialtestfn}
\end{equation}
where $NP$ is the total number of nodes in the mesh, $\bn$ is the shape function matrix and $\mathbf{d}_I$ is the vector of degrees of freedom associated with node $I$. Upon substituting \Eref{eqn:trialtestfn} into \Eref{eqn:weakform} and invoking the arbitrariness of  $\delta \uu$, we obtain the following discretized algebraic system of equations:
\begin{equation}
\kk \mathbf{d} = \mathbf{f}
\end{equation}
where
\begin{eqnarray} 
\kk &=& \int\limits_{\Omega^h} \bigb^{\rm T} \dd \bigb~\rmd \Omega \nonumber \\
\mathbf{f} &=&  \int\limits_{\Omega^h} \bn^{\rm T} \bb~\rmd \Omega + \int\limits_{\Gamma_t} \bn^{\rm T} \overline{\mathbf{t}}~\rmd \Gamma
\label{eqn:stiffmat}
\end{eqnarray}
where $\kk$ is the stiffness matrix and $\Omega^h$ is the discretized domain formed by the union of elements $\Omega^e$. The stiffness matrix is computed over each element and assembled to the global matrix. The size of the stiffness matrix depends on the number of nodes in an element.

\subsection{Generalization to arbitrary polygons}
The growing interest in the generalization of FE over arbitrary meshes has opened up a new area of finite elements called `\emph{polygonal finite elements}'. In polygonal finite elements, the number of sides of an element is not restricted to three or four as in the case of 2D. 

\begin{figure}
\centering
\includegraphics[scale=0.5]{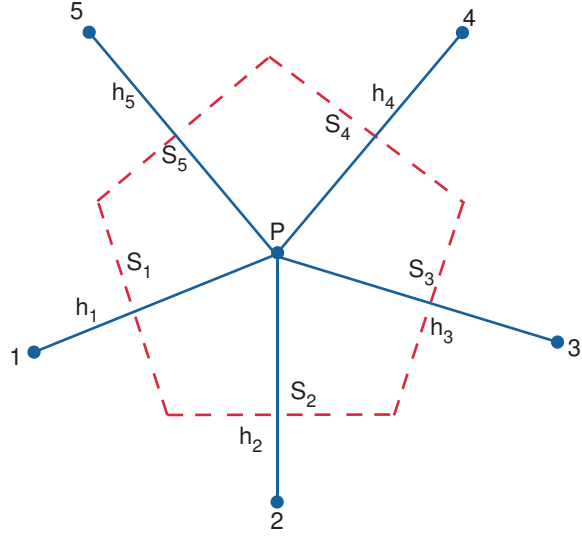}
\caption{Vorono\"{i} diagram of a point $P$.}
\label{fig:vorfig}
\end{figure}

\begin{figure}[htpb]
\centering
\subfigure[]{\includegraphics[scale=0.41]{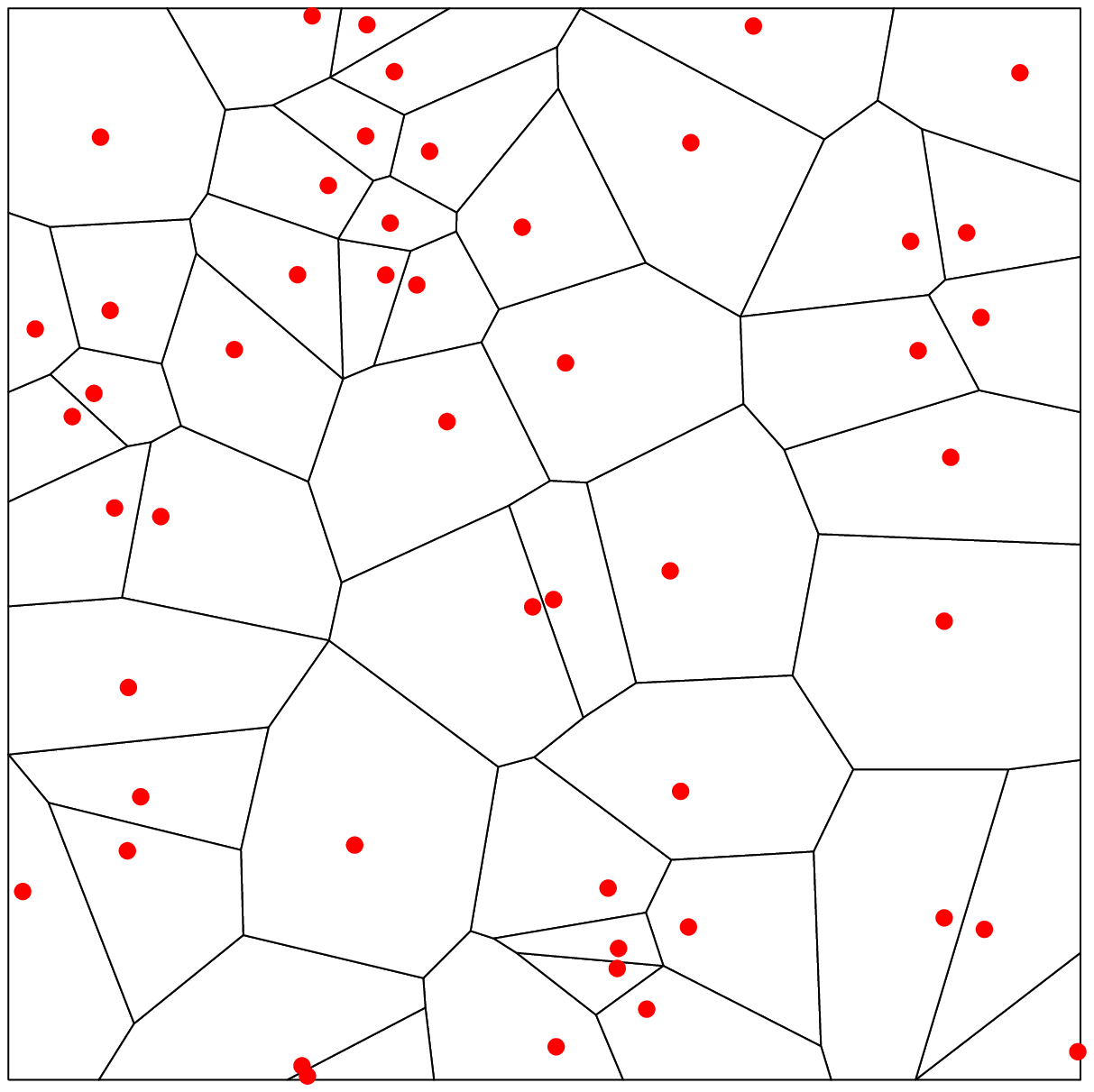}}
\subfigure[]{\includegraphics[scale=0.41]{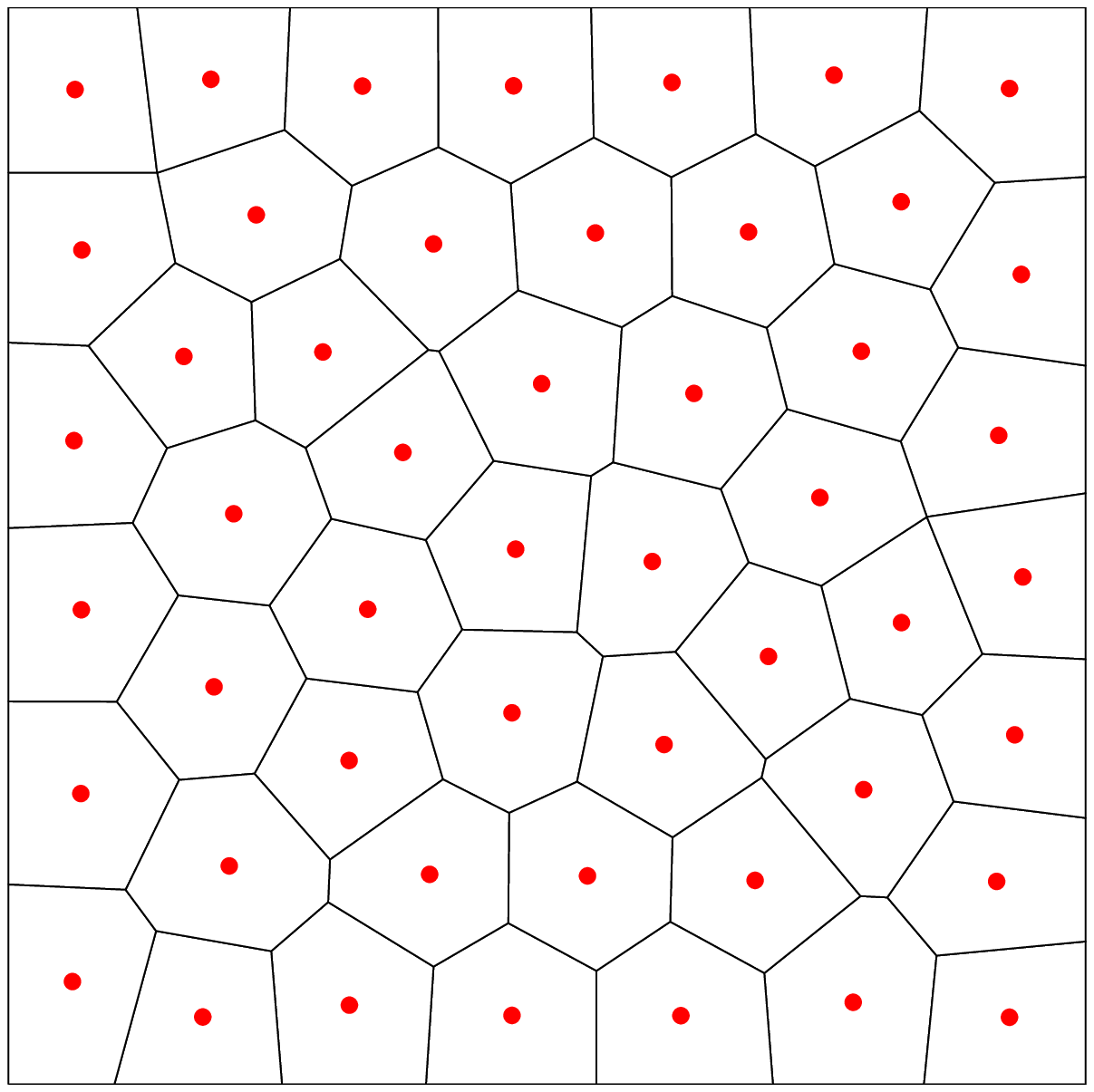}}
\caption{Vorono\"{i} tessellation of two set of data points, where the \textcolor{red}{red} dots are the seed points: (a) Scattered data set (b) Polygonal mesh after few iterations. }
\label{fig:vtessel}
\end{figure}

\begin{figure}[htpb]
\centering
\subfigure[]{\includegraphics[scale=0.25]{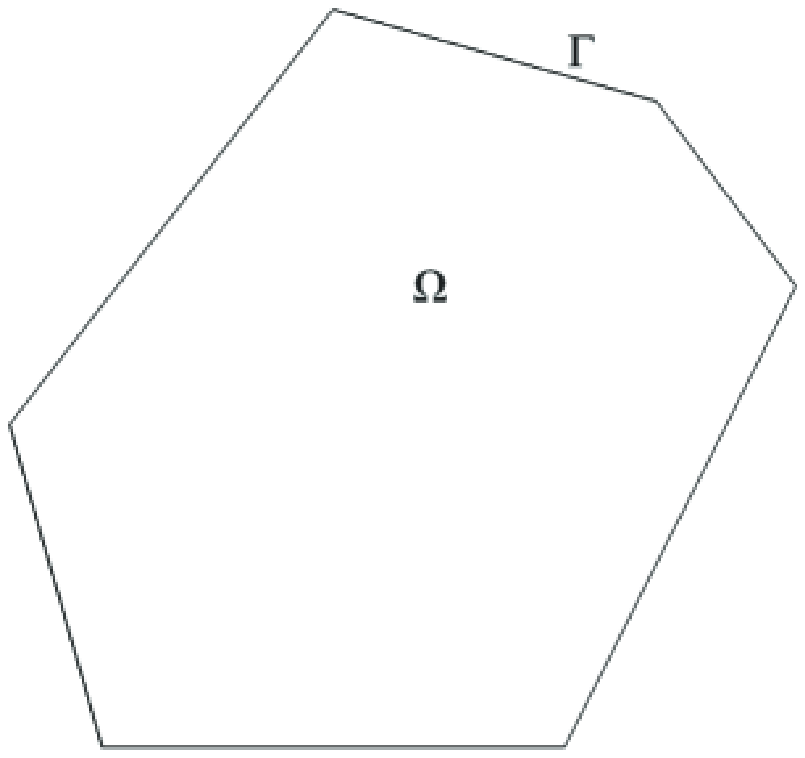}}
\subfigure[]{\includegraphics[scale=0.25]{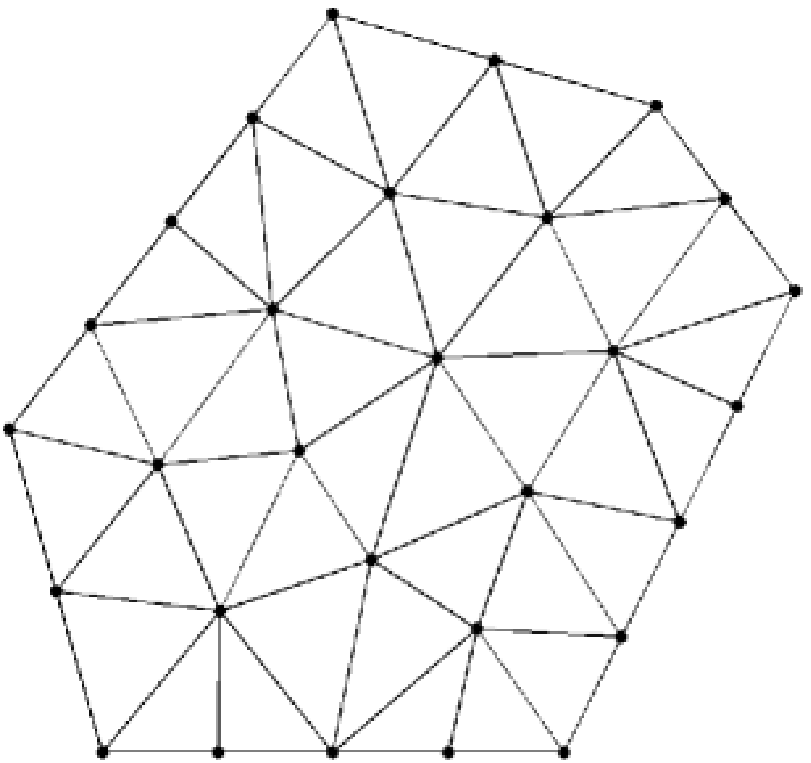}}
\subfigure[]{\includegraphics[scale=0.25]{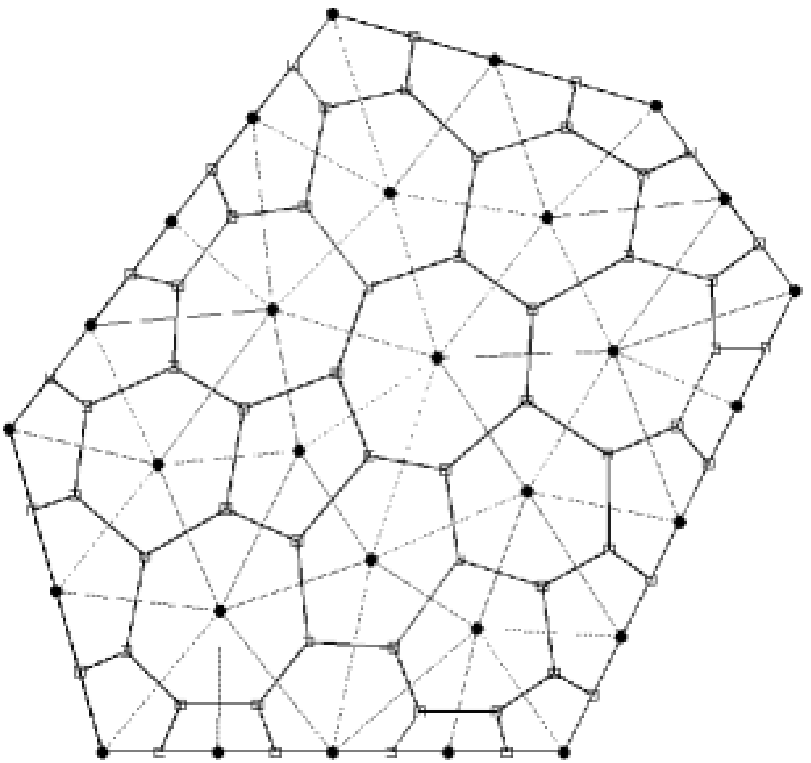}}
\caption{Polygonal mesh generation: sequence of steps.}
\label{fig:tritopoly}
\end{figure}

\paragraph{Polygonal mesh}
The Vorono\"{i} tessellation is a fundamental geometrical construct to generate a polygonal mesh covering a given domain. Polygonal meshes can be generated from Vorono\"{i} diagrams. The Vorono\"{i} diagram is a subdivision of the domain into regions $V(p_I)$, such that any point in $V(p_I)$ is closer to node $p_I$ than to any other node. \fref{fig:vorfig} shows a Vorono\"{i} diagram of a point $P$. The first order Vorono\"{i} $V(N)$ is a subdivision of the Euclidean space $\mathbb{R}^2$ into convex regions, mathematically:
\begin{equation}
T_I = \left\{ \xx \in \mathbb{R}^2 \colon d(\xx,\xx_I) < d(\xx,\xx_J) \forall J \neq I \right\}
\end{equation}
where $d(\xx_I,\xx_J)$, the Euclidean matrix, is the distance between $\xx_I$ and $\xx_J$. The quality of the generated polygonal mesh depends on the randomness in the scattered points. \fref{fig:vtessel} shows a typical Vorono\"{i} tessellation of two sets of scattered data set. The quality of a polygonal mesh determines the accuracy of the solution~\cite{sukumartabarraei2004}. To improve the quality of the Vorono\"{i} tessellation, the generating point of each Vorono\"{i} cell can be used as its center of mass, leading to a special type of Vorono\"{i} diagram, called the centroidal Vorno\"{i} tessellation (CVT)~\cite{duwang2005}. Sieger \textit{et al.,}~\cite{siegeralliez2010} presented an optimizing technique to improve the Vorono\"{i} diagrams for use in FE computations. A polygonal mesh can also be constructed from a triangular mesh by connecting the centroid of  all the triangular elements circumventing a particular node~\cite{ooisong2012}. \fref{fig:tritopoly} shows a sequence of steps required to generate a polygonal mesh from an underlying triangular mesh. This procedure is also adopted in the edge based smoothed FEM~\cite{nguyen-thoiliu2011}. In the following sub-sections we briefly discuss three different displacement based formulations.

\subsection{Polygonal Finite Element Method (Polygonal FEM)}
In polygonal finite elements, the domain is discretized with elements having arbitrary number of edges. The shape functions are constructed over each polygonal element. The construction of shape functions over arbitrary polygons can be broadly categorized as:
\begin{itemize}
\item Length and area measures (Wachspress interpolants)~\cite{wachspress1971,floatergillette2013}.
\item Natural neighbour interpolants (Sibson, Laplace interpolant)~\cite{sibson1980,sukumarmoran1998}.
\item Harmonic coordinates~\cite{warrenschaefer2007,bishop2013}.
\item Maximum entropy approximants~\cite{sukumar2013}.
\item Mean value coordinates~\cite{floater2003}.
\end{itemize}

For recent advances in polygonal elements, interested readers are referred to the literature~\cite{sukumarmalsch2006} and the references therein. In this study, we consider only Laplace interpolants with isoparametric mapping\footnote{When using the isoparametric formulation, the Wachspress interpolants and the Laplace interpolants are identical~\cite{sukumartabarraei2004}. Hence, the interpretation of the results with the Laplace interpolants are applicable when Wachspress interpolants are used.}. The Laplace interpolant is also called the natural neighbor interpolant~\cite{sibson1980,sukumarmoran1998,sukumartabarraei2004}. It provides a natural weighting function for irregularly spaced nodes. For a point $P$ with $n$ natural neighbors, the Laplace shape functions for node $P_I$ can be written as~\cite{sukumartabarraei2004}:
\begin{equation}
\phi_I(\xx) = \frac{ \alpha_I(\xx)}{\sum\limits_{I=1}^n \alpha_J(\xx)}, \hspace{0.25cm} \alpha_J(\xx) = \frac{s_J(\xx)}{h_J(\xx)}, \hspace{0.2cm} \xx \in \mathbb{R}^2
\end{equation}
where $\alpha_I(\xx)$ is the Laplace weight function, $s_I(\xx)$ is the length of the Vorono\"{i} edge associated with $P$ and node $P_I$ and $h_I(\xx)$ is the Euclidean distance between $P$ and $P_I$ (see \fref{fig:vorfig}). Listing \ref{laplaceshapefn} gives a MATLAB \textsuperscript{\textregistered} function to compute the Laplace interpolants for a point $P$\footnote{Original source code in FORTRAN is available from~\cite{sukumarlaplace}}.

\begin{lstlisting}[caption=Laplace Interpolant over convex polytope, label=laplaceshapefn]
function [phi,dphix,dphiy]=shape(x,y,coord)
shapedenom=0; n=length(coord(:,1)); % number of sides of the polygon
phi = zeros(n,1); sumder = zeros(n,1);
for i=1:n %loop over number of vertices
    a1 = coord(i,1);  a2 = coord(i,2);
    if(i==1)
        nodem = n;  nodep = i+1;
    elseif(i==n)
        nodem = n-1; nodep = 1;
    else
        nodem = i-1;  nodep = i+1;
    end
    %circumcenter of a triangle formed by points: (a1,b1,x) & (a2,b2,y)
    b1 = coord(nodem,1);  b2 = coord(nodem,2);    
    [v1,v2,v1x,v1y,v2x,v2y] = circum(b1,b2,a1,a2,x,y);
    %circumcenter of a triangle formed by points: (a1,c1,x) & (a2,c2,y)
    c1 = coord(nodep,1);  c2 = coord(nodep,2);
    [s1,s2,s1x,s1y,s2x,s2y]=circum(a1,a2,c1,c2,x,y);
    SI = sqrt((v1-s1)^2+(v2-s2)^2); HI = sqrt((x-a1)^2+(y-a2)^2);
    SI_v1 = ((v1-s1)*(v1x-s1x) + (v2-s2)*(v2x - s2x))/SI;
    SI_v2 = ((v1-s1)*(v1y-s1y) + (v2-s2)*(v2y - s2y))/SI; 
    HI_v1 = (x-a1)/HI; HI_v2 = (y-a2)/HI;
    shapenum = SI/HI;
    der(1) = (SI_v1 - shapenum*HI_v1)/HI;
    der(2) = (SI_v2 - shapenum*HI_v2)/HI;  
    shapedenom = shapedenom + shapenum;
    sumder(1) = sumder(1) + der(1); sumder(2) = sumder(2) + der(2);   
    shapefun(i) = shapenum;
    phix(i) = der(1); phiy(i) = der(2);
end

phi = shapefun./shapedenom;
dphix = (phix - phi.*sumder(1))/shapedenom;
dphiy = (phiy - phi.*sumder(2))/shapedenom;

function [v1,v2,v1x,v1y,v2x,v2y] = circum(a1,a2,b1,b2,x,y)
D = (a1-x)*(b2-y) - (b1-x)*(a2-y);
t1 = ((a1-x)*(a1+x) + (a2-y)*(a2+y))*(b2-y);
t2 = ((b1-x)*(b1+x) + (b2-y)*(b2+y))*(a2-y);
term1 = (1/2)*(t1-t2); v1 = term1/D;
t4 = ((b1-x)*(b1+x) + (b2-y)*(b2+y))*(a1-x);
t5 = ((a1-x)*(a1+x) + (a2-y)*(a2+y))*(b1-x);
term2 = (1/2)*(t4-t5); v2 = term2/D;
dx = a2-b2; dy = b1-a1;
term3 = (1/2)*((b1^2-a1^2)+(b2^2-a2^2));
v1x = (x-v1)*dx/D;
v1y = (term3 + y*dx - v1*dy)/D;
v2x = (-term3 + x*dy - v2*dx)/D;
v2y = (y-v2)*dy/D;
end
\end{lstlisting}

\subsubsection{Numerical integration} The computation of the stiffness matrix and the force vectors involves evaluating the integrals given by \Eref{eqn:stiffmat}. The shape functions over arbitrary polygons are rational polynomials and hence the formulation of efficient integration rules poses a unique challenge. This has received much attention in the recent years~\cite{sukumartabarraei2004,natarajanbordas2009,mousavixiao2010,nguyen-thoiliu2011} as evidenced by the growing literature. One potential solution is to sub-divide the physical element into triangles and then use the well-known quadrature on a triangle for numerical integration~\cite{sukumartabarraei2004}. The purpose of subdivision is solely for the purpose of numerical integration and does not introduce additional unknowns. Although, straightforward, this process involves a two-level isoparameteric mapping and relies on the positivity of the Jacobian matrix involved in the transformation. Lyness and Monegato~\cite{lynessmonegato1977} presented quadrature rules for regions with hexagonal symmetry. Natarajan~\textit{et al.,}~\cite{natarajanbordas2009} proposed a numerical integration technique over arbitrary polygons based on complex mapping. This procedure eliminates the need for a two level isoparametric mapping while guaranteeing the positivity of the Jacobian, but is restricted to only 2D. Sommariva and Vianello~\cite{sommarivavianello2009} presented a Gauss-like cubature over arbitrary polygons. Mousavi \textit{et al.,}~\cite{mousavixiao2010} presented a numerical algorithm based on group theory and optimization scheme to compute quadrature rules over arbitrary polygons. 

Once the stiffness matrix and the force vector over each element is computed, then the conventional FE procedure is adopted for the assemblage process and the solution of the system of equations~\cite{sukumartabarraei2004,natarajanbordas2009}. Listing \ref{pfemstiffmat} presents the MATLAB \textsuperscript{\textregistered} function to compute the stiffness matrix. The shape functions are computed by Laplace interpolants using the routine given in Listing \ref{laplaceshapefn}.

\begin{lstlisting}[caption=Computation of the stiffness matrix using conventional Polygonal FEM, label=pfemstiffmat]
function [kmat] = getKmatPolyFEM(element,node,matmtx,W,Q)
for iel = 1:size(element,1)    
    econ = element{iel}; nn=length(econ); ndof=2;
    gindex = []; gindex=reshape([ndof*econ-1;ndof*econ],1,[]);
    for igp = 1:size(W,1)      % loop over integration points
    	% shape function and derivatives in canonical domain
        [N,dphir,dphis] = shape(Q(igp,1),Q(igp,2),polycoord(length(econ)) ;) ;  
        JO = [dphir'; dphis']*node(econ,:) ;
        dNdx = inv(JO)*[dphir'; dphis'] ;
        % strain-displacement matrix
        B = zeros(3,ndof*nn) ;  dNdx = dNdx' ;
        B(1,1:ndof:ndof*nn) = dNdx(:,1)' ; B(2,2:ndof:ndof*nn) = dNdx(:,2)' ;
        B(3,1:ndof:ndof*nn) = dNdx(:,2)' ; B(3,2:ndof:ndof*nn) = dNdx(:,1)' ;
        % compute the stiffness matrix
        kmat(gindex,gindex) = kmat(gindex,gindex) + B'*matmtx*B*W(igp)*det(JO) ;
    end  % end loop over integration points
end  % end loop over elements
\end{lstlisting}

\subsection{Polygonal Smoothed Finite Element Method (Polygonal SFEM) }
Similar to the Polygonal FEM, the Polygonal SFEM, also discretizes the domain into polygonal elements with arbitrary edges. However, the stiffness matrix given by \Eref{eqn:stiffmat}, is computed by using the strain projection procedure. In this approach, the strain field, $\tilde{\varepsilon}_{ij}^h$ used to compute the stiffness matrix is computed from a weighted average of the standard strain field $\varepsilon_{ij}^h$. At a point $\xx_{C}$ in an element $\Omega^h$, the smoothed strain field is given by:
\begin{align}
\tilde{\varepsilon}_{ij}^h &= \int\limits_{\Omega^h} \varepsilon_{ij}^h (\xx) \Phi(\xx - \xx_C)~\mathrm{d}\xx \nonumber \\
&= \int\limits_{\Omega_C} \bigb \qq \Phi(\xx - \xx_C)~\mathrm{d}\xx = \widetilde{\bigb} \qq
\end{align}
where $\qq$ is a vector of degrees of freedom, $\Phi$ is a smoothing function is given by
\begin{equation*}
\renewcommand{\arraystretch}{1.0}
\Phi = \left\{ \begin{array}{cc} \frac{1}{A_C} & \xx_C \in \Omega_C \\ 0 & \xx_C \notin \Omega_C \end{array} \right.
\end{equation*}
where $A_C$ is the area of the subcell. The smoothed strain-displacement matrix $\tilde{\bigb}$ is given by:
\begin{equation}
\renewcommand{\arraystretch}{1.5}
\widetilde{\bigb} = \frac{1}{A_C} \int\limits_{\Omega_C} \left[ \begin{array}{cc} \frac{\partial N_I}{\partial x} & 0 \\ 0 & \frac{\partial N_I}{\partial y} \\ \frac{\partial N_I}{\partial y} & \frac{\partial N_I}{\partial x} \end{array} \right] ~\rmd \Omega
\end{equation}
Using the divergence theorem, we obtain
\begin{equation}
\renewcommand{\arraystretch}{1.5}
\widetilde{\bigb} = \frac{1}{A_C} \int\limits_{\Gamma_C} \left[ \begin{array}{cc} n_x N_I & 0 \\ 0 & n_y N_I \\ n_y N_I & n_x N_I \end{array} \right] ~\rmd \Gamma
\label{eqn:smoothedB}
\end{equation}
where $N_I$ is the shape function at node $I$. The smoothed stiffness matrix is given by~\cite{liutrung2010}:
\begin{equation}
\widetilde{\kk} = \widetilde{\bigb}^{\rm T} \dd \widetilde{\bigb} A_C
\end{equation}
This process of smoothing the gradient field is known as `\emph{cell-based smoothed finite element method}' (CSFEM) in the literature~\cite{liunguyen2007}. Until now, no assumption has been made on the shape of the element. The procedure outlined so far is general and is applicable to polygons of arbitrary shapes~\cite{dailiu2007}. Other smoothing techniques~\cite{liutrung2010} can also be extended to polygons, for example, Nguyen-Thoi \textit{et al.,}~\cite{nguyen-thoiliu2011} applied the edge based smoothing technique over polygons for problems in solid mechanics. In this study we study the properties of the CSFEM over arbitrary polygons.

\begin{figure}[htpb]
\centering
\includegraphics[scale=0.75]{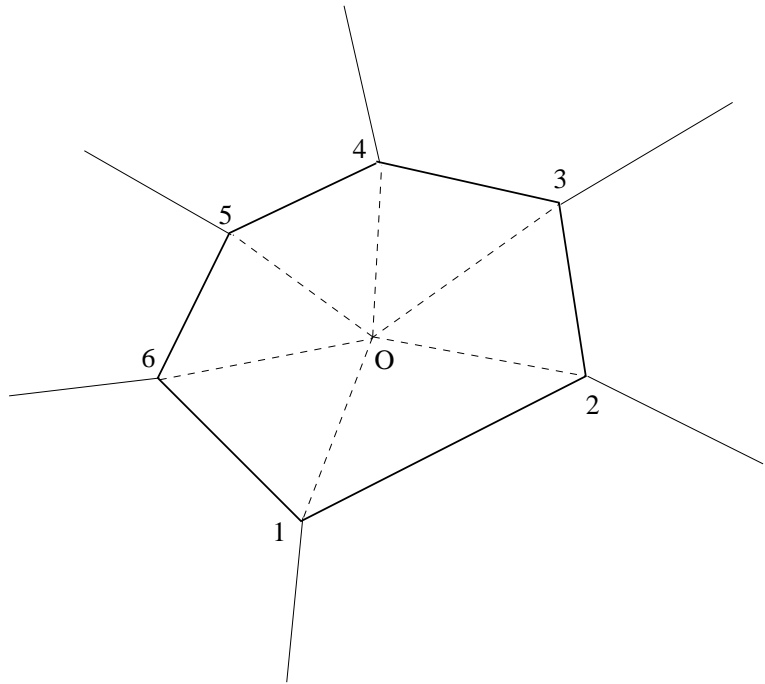}
\caption{Construction of simple averaging shape functions for a polygonal element. $O-1-2$ is subcell.}
\label{fig:smoothedavgshpfn}
\end{figure}

\paragraph{Construction of Shape functions} Due to the process of strain smoothing, only the shape function is involved in the calculation of the field gradients and hence the stiffness matrix. Moreover, no analytical form of the shape functions or its derivatives is required. In this study, we use the simple averaging technique outlined in~\cite{dailiu2007} to compute the shape function over arbitrary polygons. The construction of the shape functions is as follows: for a general polygonal element, a central point $O$ is located by (see \fref{fig:smoothedavgshpfn}):
\begin{equation}
(x_O,y_O) = \frac{1}{n} \sum\limits_{I}^n(x_I,y_I)
\end{equation}
where $n$ is the number of nodes of the polygonal element. The shape function at point $O$ is given by $\left[ \frac{1}{n} \cdots \frac{1}{n} \right]$ with size 1$\times n$. In this study, the shape function along each edge is linear and one Gauss point is sufficient to evaluate the integral in \Eref{eqn:smoothedB}. For the midpoint on the side of the element, its shape functions are evaluated averagely using two related nodes while for the interior point, its shape functions are evaluated using the point $O$ and the related node.
Once the shape functions are evaluated, the stiffness matrix can be computed by following the procedure outlined in~\cite{liunguyen2007,dailiu2007}.

\paragraph{Choice of number of subcells} When the strain smoothing is applied to quadrilateral elements, each element can be subdivided into number of triangular or quadrilateral subcells. Liu~\textit{et al.,}~\cite{liunguyen2007} showed that upper and lower bound solutions in SFEM can be obtained by varying the number of smoothing cells. In case of arbitrary polygons, Dai \textit{et al.,}~\cite{dailiu2007} showed that for the stiffness matrix to be non-singular and stable, an $n-$sided polygon should be sub-divided into $n$ triangular subcells. In this study, we adopt $n$ triangular subcells for the construction of the stiffness matrix. A simple Matlab\textsuperscript{\textregistered} function to compute the smoothed stiffness matrix is given in Listing \ref{nsfemstiffmat}.

\begin{lstlisting}[caption=Smoothed stiffness matrix, label=nsfemstiffmat]
function [kmat] = getKmatPolySFEM(element,node,matmtx,W,Q)
for iel = 1:numelem
    gecon = element{iel} ; % get current element connectivity
    coord = node(gecon,:); nn = lenght(econ); ndof = 2;
    tri = delaunay(ncoord(:,1),ncoord(:,2));  % triangulate to get subcells
    N = [eye(nn); 1/nn*ones(1,nn)]; % shape function at nodes and at the centroid
    gindex = []; gindex=reshape([ndof*econ-1;ndof*econ],1,[]);
    for isc = 1:size(tri,1)   % loop over subcells...
        B = zeros(nn,2);
        ncon = tri(isc,:); nodes = length(ncon) ;
        xl = xd(ncon,1); yl = xd(ncon,2);
        [geom,iner,cpmo] = polygeom(xl,yl); % compute the area of the subcell...
        for iside = 1:size(tri,2)   % loop over the sides of the subcell
            xy1 = [xl; xl(1)]; xy2 = [yl; yl(1)]; xys = [xy1 xy2];
            shape = [N(ncon,:); N(ncon(1),:)];  % shape function 
            lside = norm( xys(iside+1,:)-xys(isde,:) );   % side length
            nx = 1/2*(xy2(iside+1)-xy2(iside))/(lside/2) ;  % normals
            ny = -1/2*(xy1(iside+1)-xy1(iside))/(lside/2) ;
            phi = mean( shape(iside,:),shape(iside+1,:));  % shape function            
            % construct the B matrix...
            for inode = 1:nn
                B(inode,1) = B(inode,1) + phi(inode)*nx*lside/geom(1);
                B(inode,2) = B(inode,2) + phi(inode)*ny*lside/geom(1);
            end
        end
        Bmat(1:ndof:ndof*nn,1) = B(:,1); Bmat(2:ndof:ndof*nn,2) = B(:,1);
        Bmat(1:ndof:ndof*nn,3) = B(:,2); Bmat(2:ndof:ndof*nn,3) = B(:,1);
        % stiffness matrix    
        kmat(gindex,gindex) = kmat(gindex,gindex) +Bmat*matmtx*Bmat'*geom(1) ;
    end
end
\end{lstlisting}

\subsection{Polygonal Scaled Boundary Finite Element Method (Polygonal SBFEM)}
Polygonal elements can be conveniently developed from the scaled boundary FEM (SBFEM). The SBFEM~\cite{wolfsong2001} is a semi-analytical computational technique that reduces the governing partial differential equations to a set of ordinary differential equations. In the SBFEM, a local radial-circumferential-like coordinate system is introduced. Numerical solutions are sought around the circumferential direction using conventional FEM, whilst in the radial direction, the solution is defined by smooth analytical functions. A scaling centre $O$ is selected at a point from which the whole boundary of the domain is visible. This `\emph{scaling requirement}' can always be satisfied by sub-structuring, i.e. dividing the structure into smaller subdomains. 

\begin{figure}
\centering
\subfigure[Uncracked polygon]{\includegraphics[height=7cm]{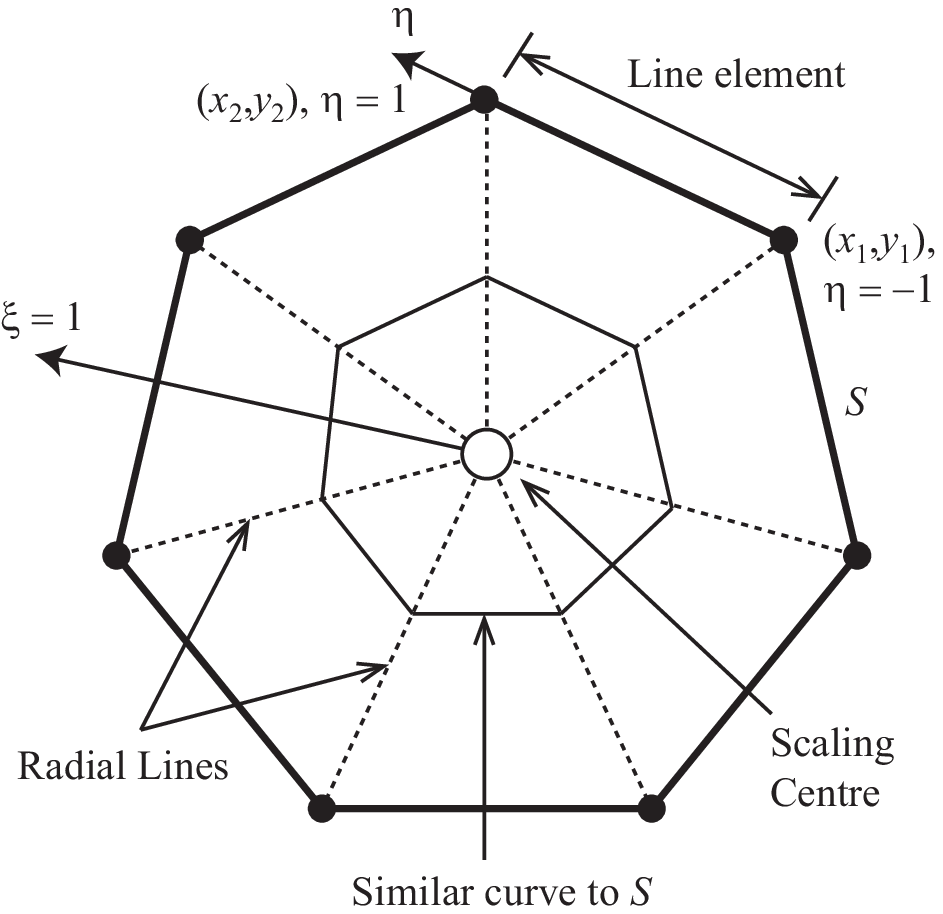}}
\subfigure[Cracked polygon]{\includegraphics[height=7cm]{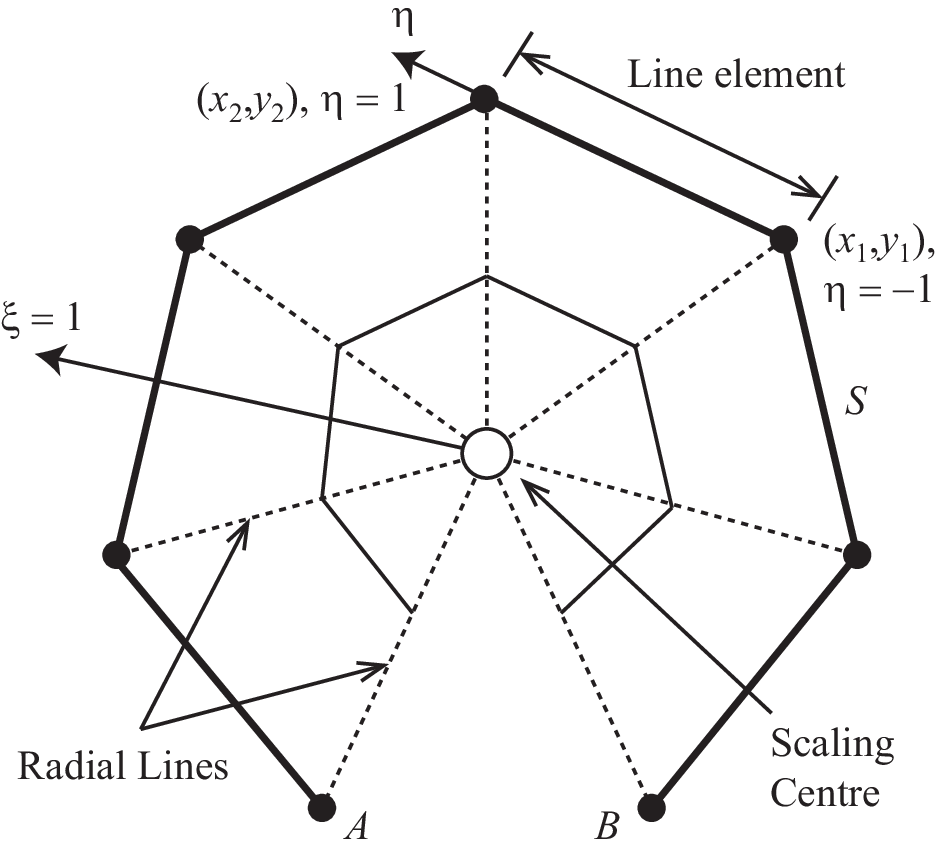}}
\caption{Polygon representation by the SBFEM.}
\label{fig:Polygon-representation-by}
\end{figure}

\fref{fig:Polygon-representation-by} shows a polygon modeled by the SBFEM. There is no restriction to the number of sides that a polygon can have. The radial coordinate $\xi$ is defined from the scaling centre $(\xi=0)$ to the polygon boundary $(\xi=1)$. Each edge on the polygon boundary is discretised using 1D finite elements with local coordinate $\eta$. The local coordinate is defined in the range $-1\leq\eta\leq1$ and is similar to the standard finite elements.

\begin{remark} The geometry of the polygon has to only satisfy the aforementioned scaling requirement. This condition is automatically satisfied for all convex polygons and many concave polygons. 
\end{remark}

\paragraph{Displacement approximation} The displacements of a point in a polygon is approximated by: 
\begin{equation}
\mathbf{u}(\xi,\eta)=\mathbf{N}(\eta)\mathbf{u}(\xi)\label{eqn:dispapprox}
\end{equation}
where $\mathbf{u}(\xi)$ are radial displacement functions and $\mathbf{N}(\eta)$ is the shape function matrix in the circumferential direction. Along the circumferential direction, the standard 1D Gauss-Lobatto-Lagrange shape functions can be used. 
\begin{remark}The shape functions $\mathbf{N}(\eta)$ can be of any order and does not affect the way in which the stiffness matrix and load vectors are evaluated. Higher order polygons to be constructed conveniently using the SBFEM by simply increasing the order of the shape functions in $\mathbf{N}(\eta)$.
\end{remark}
 
By following the procedure outlined in~\cite{wolfsong2001,deekswolf2002}, the following ODE is obtained:
\begin{equation}
\mathbf{E}_{0}\xi^{2}\mathbf{u}(\xi)_{,\xi\xi}+(\mathbf{E}_{0}+\mathbf{E}_{1}^{\mathrm{T}}-\mathbf{E}_{1})\xi\mathbf{u}(\xi)_{,\xi}-\mathbf{E}_{2}\mathbf{u}(\xi)=0\label{eqn:governODEsbfem}
\end{equation}
where $\mathbf{E}_{0},\mathbf{E}_{1}$ and $\mathbf{E}_{2}$ are coefficient matrices given by:
\begin{align}
\mathbf{E}_{0} & =\int_{\eta}\mathbf{B}_{1}(\eta)^{{\rm T}}\mathbf{D}\mathbf{B}_{1}(\eta)|\mathbf{J}(\eta)|d\eta,\nonumber \\
\mathbf{E}_{1} & =\int_{\eta}\mathbf{B}_{2}(\eta)^{{\rm T}}\mathbf{D}\mathbf{B}_{1}(\eta)|\mathbf{J}(\eta)|d\eta,\nonumber \\
\mathbf{E}_{2} & =\int_{\eta}\mathbf{B}_{2}(\eta)^{{\rm T}}\mathbf{D}\mathbf{B}_{2}(\eta)|\mathbf{J}(\eta)|d\eta.\label{eqn:coeffmat}
\end{align}
They are evaluated element-by-element on the polygon boundary and assembled over a polygon. This process is similar to the standard FE procedure of assemblage. The boundary nodal forces are related to the displacement functions by: 
\begin{equation}
\mathbf{f}=(\mathbf{E}_{0}\xi\mathbf{u}(\xi)_{,\xi}+\mathbf{E}_{1}^{{\rm T}}\mathbf{u}(\xi))|_{\xi=1}\label{eqn:nodalforce}
\end{equation}
The stiffness matrix $\mathbf{K}$ in case of the SBFEM is given by:
\begin{equation}
\mathbf{K}=\boldsymbol{\Phi}_{\mathrm{q}}\boldsymbol{\Phi}_{\mathrm{u}}^{-1}\label{eqn:sbfemkmat-b}
\end{equation}
where $\boldsymbol{\Phi}_{\mathrm{u}}$ and $\boldsymbol{\Phi}_{\mathrm{q}}$ represent the modal displacements and forces, respectively, computed from the eigenvalue decomposition of the following Hamiltonian matrix~\cite{wolfsong2001}:
\begin{equation}
\mathbf{Z}=\left[\begin{array}{cc}
\mathbf{E}_{0}^{-1}\mathbf{E}_{1}^{\mathrm{T}} & -\mathbf{E}_{0}^{-1}\\
\mathbf{E}_{1}\mathbf{E}_{0}^{-1}\mathbf{E}_{1}^{\mathrm{T}}-\mathbf{E}_{2} & -\mathbf{E}_{1}\mathbf{E}_{0}^{-1}
\end{array}\right]\label{eq:Hamiltonian matrix}
\end{equation}
The corresponding eigenvalues are given by: $\boldsymbol{\Lambda}_{\mathrm{n}}=\mathrm{diag}\left(\lambda_{1},\,\lambda_{2},\,...,\lambda_{n}\right)$. For a 2 noded line element, the coefficient matrices given by~\Eref{eqn:coeffmat} can be written explicitly. The one dimensional shape functions up to order 3 is given in Listing \ref{shapefn}. A simple Matlab\textsuperscript{\textregistered} function is given in Listing \ref{sbfemstiffmat} to compute the coefficient matrices and the corresponding stiffness matrix.

\begin{lstlisting}[caption=One dimensional shape functions, label=shapefn]
function [shp] = shapeFunction1d(ennodes, x)
%  shp(1,j):  shape function Nj(eta) of j-th node
%  shp(2,j):  derivative of shape function Nj(eta),eta
%  ennodes: number of nodes of element. x:        local coordinate eta
shp = zeros(2,ennodes);
switch ennodes
    case(2) % 2-node line element
       shp = [0.5d0*(1.d0-x) 0.5d0*(1.d0+x);-0.5d0 0.5d0];        
    case(3) %  3-node line element
        shp = [0.5d0*(-1.d0+x)*x 1.d0-x*x 0.5d0*(1.d0+x)*x;-0.5d0+x -2*x 0.5d0+x];       
    case(4)
        % Shape functions of 4-node line element
        shp(1,1)=-0.125d0*(-1.d0+x)*(-1.d0+5.d0*x*x);
        shp(1,2)= 0.625d0*(-1.d0+x*x)*(-1.d0+2.2360679774997897d0*x);
        shp(1,3)=-0.27950849718747371d0*(-1.d0+x*x)*(2.2360679774997897d0+5.d0*x);
        shp(1,4)= 0.125d0*( 1.d0+x)*(-1.d0+5.d0*x*x);
        shp(2,1)= 0.125d0*( 1.d0+(10.d0-15.d0*x)*x);
        shp(2,2)= 0.625d0*(-2.2360679774997897d0+(-2.d0+6.7082039324993691d0*x)*x);
        shp(2,3)=-0.279508497187473712d0*(-5.d0+(4.472135954999579393d0+15.d0*x)*x);
        shp(2,4)= 0.125d0*(-1.d0+(10.d0+15.d0*x)*x);
end
\end{lstlisting}

\begin{remark} When the SBFEM is applied to arbitrary polygons, the stiffness matrix is computed directly using \Eref{eqn:sbfemkmat-b}. An explicit form of the shape functions or its derivatives are not required for linear problems. For nonlinear problems, the shape functions can be derived by solving the the elasticity equation~\cite{ooisong2013}.
\end{remark}

\begin{lstlisting}[caption=Scaled boundary polygon formulation: stiffness matrix, label=sbfemstiffmat]
function [kmat] = getKmatPolySBFEM(element,node,matmtx,W,Q)
for iel = 1:numelem	   % loop over the elements    
    econ = element{iel}; nn=length(econ); ndof=2;
    gindex = []; gindex=reshape([ndof*econ-1;ndof*econ],1,[]);
    EleUnknown = length(econ)*ndof; E0 = zeros(EleUnknown,EleUnknown) ;
    E1 = zeros(EleUnknown,EleUnknown) ; E2 = zeros(EleUnknown,EleUnknown) ;
    for isub = 1:size(subecon,1)  % loop over the sides of the polygon
        lecon = subecon(isub,:) ; xynodes = xyc(lecon,:) ;
        lindex = reshape([ndof*lecon-1;ndof*lecon],1,[]);
        for igp = 1:size(W,1)  % loop over the integration points...
            [shp] = shapeFunction1d(length(lecon), Q(igp,:)) ;  
            N = shp(1,:); dNdxi = shp(2,:);  % shape functions and its derivatives
            Jac = [N' dNdxi']'*xynodes(:,1:2) ;		% Jacobian 
            % small b1 and b2 terms...
            xetaN = Jac(2,1)*N; yetaN = Jac(2,2)*N;
            xetadN = Jac(1,1)*dNdxi; yetadN = Jac(1,2)*dNdxi;
            
            b1eta = [yetaN; zeros(1,length(N)); -xetaN; zeros(1,length(N));-xetaN; yetaN];
            b1eta = 1/det(Jac)*reshape(b1eta,3,[]);
            b2eta = [-yetadN; zeros(1,length(N)); xetadN; zeros(1,length(N));xetadN; -yetadN];
            b2eta = 1/det(Jac)*reshape(b2eta,3,[]);
            
            % coefficient matrices.... lindex - local index
            E0(lindex,lindex)=E0(lindex,lindex)+b1eta'*matmtx*b1eta*det(Jac)*W(igp);            
            E1(lindex,lindex)=E1(lindex,lindex)+b2eta'*matmtx*b1eta*det(Jac)*W(igp);            
            E2(lindex,lindex)=E2(lindex,lindex)+b2eta'*matmtx*b2eta*det(Jac)*W(igp);
        end
    end
    % compute the stiffness matrix...
    nd = ndof*nn; id = 1:nd;
    m = E0\[E1' -eye(nd)]; Z = [m; E1*m(:,id)-E2 E1*m(:,nd+id)];
    [v, d] = eig(Z);  lambda = diag(d); [~, idx] = sort(real(lambda),'ascend'); 
    lambda = lambda(idx(id));  v = v(:, idx(id)); 
    v(:,end-1:end) = 0; v(1:2:nd,  end-1) = 1; v((1:2:nd)+1,end) = 1;
    % stiffness matrix...
    kmat(gindex,gindex) = kmat(gindex,gindex) + real(v(nd+id, :)/v(id, :)) ;
end
\end{lstlisting}

\section{Numerical Examples}
\label{numexamp}
In the first part of this section, the results from the three different finite element techniques discussed above are compared using three benchmark problems in the context of linear elasticity. The results are compared with the analytical solution where available. In the later part of the section, problems involving strong discontinuity (e.g. plate with a crack, two cracks emanating from a hole) are solved using the scaled boundary polygon formulation. The results from the scaled boundary polygon formulation are compared with the results from other techniques taken from literature~\cite{tabarraeisukumar2008,saimotomotomura2010,dauxmoes2000}. We employ the following convention while discussing the results:
\begin{itemize}
\item Polygonal FEM - conventional polygonal FEM with Laplace interpolants. For the purpose of numerical integration, sub-triangulation is employed and sixth order Dunavant quadrature rule~\cite{dunavant1985} over each triangle is employed.
\item Polygonal nSFEM - $n-$sided smoothed finite element method. A simple averaging technique is employed to compute the shape functions and the corresponding stiffness matrix. In this study, along each edge, the shape functions are assumed to be linear and hence, only one Gauss point is sufficient to integrate the terms in the stiffness matrix. 
\item Polygonal SBFEM - scaled boundary polygon fomulation is employed over each polygon. In this case, only the polygon boundary is discretized and 1D shape functions along the boundary are used to approximate the displacement field. Hence, no special numerical integration technique is required to compute the stiffness matrix.
\end{itemize}

The built-in Matlab\textsuperscript{\textregistered} function {\small voronoin} and Matlab\textsuperscript{\textregistered} functions in {\small PolyTop}~\cite{talischipaulino2012} for building the mesh-connectivity are used to create the polygonal meshes. For the purpose of error estimation and convergence studies, the relative error, $L^2$ and $H^1$ norms, are used. The displacement norm is given by:

\begin{equation}
|| \uu - \uu^h||_{L^2(\Omega)} = \sqrt{ \int\limits_{\Omega} [ (\uu - \uu^h) \cdot (\uu-\uu^h)]~\rmd \Omega}
\end{equation}
where $\uu^h$ is the numerical solution and $\uu$ is the analytical solution or a reference solution. The energy norm is given by:

\begin{equation}
|| \uu - \uu^h||_{H^1(\Omega)} = \sqrt{ \int\limits_{\Omega} [ (\bveps - \bveps^h)^{\rm T} \dd (\bveps-\bveps^h)]~\rmd \Omega}
\end{equation}

\subsection{Patch test} 
\begin{figure}[htpb]
\centering
\scalebox{0.80}{\input{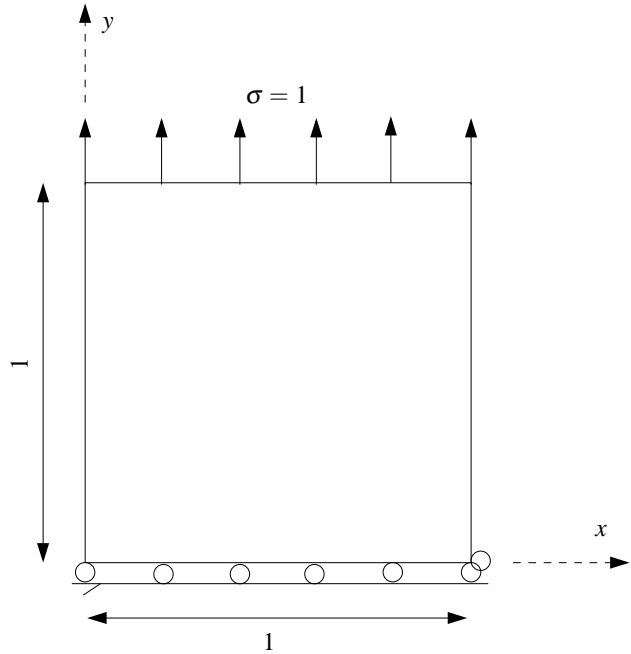}}
\caption{Geometry and loads for equilibrium patch test.}
\label{fig:simplepatch}
\end{figure}

In this section, the equilibrium path test is performed to study the accuracy and the convergence properties of different finite element techniques. Consider a uniaxial stress $\sigma=$ 1, under the condition of plate stress applied in the $y-$direction on the top edge of a unit square, whilst essential boundary conditions are applied at the bottom edge as shown in \fref{fig:simplepatch}. The exact displacement field is given by:

\begin{align}
u(x,y) &= \frac{\nu}{E}(1-x) \nonumber \\
v(x,y) &=\frac{y}{E}
\end{align}

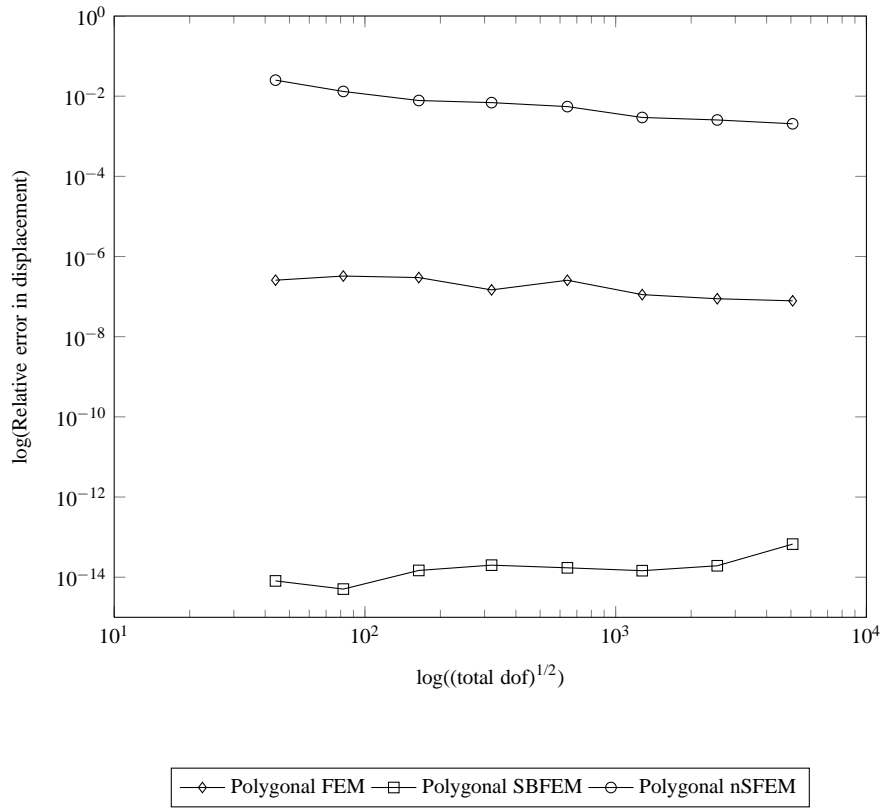
\begin{figure}
\centering
\newlength\figureheight 
\newlength\figurewidth 
\setlength\figureheight{8cm} 
\setlength\figurewidth{10cm}
%
%
%
%
\begin{tikzpicture}

\begin{axis}[%
width=\figurewidth,
height=\figureheight,
scale only axis,
xmode=log,
xmin=10,
xmax=10000,
xminorticks=true,
xlabel={$\text{log((total dof)}^{\text{1/2}}\text{)}$},
ymode=log,
ymin=1e-15,
ymax=1,
yminorticks=true,
ylabel={log(Relative error in displacement)},
legend style={at={(0.5,-0.25)},
      anchor=north,legend columns=-1},
]
\addplot [
color=black,
solid,
mark=diamond,
mark options={solid}
]
table[row sep=crcr]{
44 2.56911438424756e-07\\
82 3.25824717298471e-07\\
164 2.98225824126995e-07\\
320 1.47057104359798e-07\\
642 2.55924435178871e-07\\
1276 1.12040981943283e-07\\
2542 8.83744414918274e-08\\
5080 7.84573770824522e-08\\
};
\addlegendentry{Polygonal FEM};

\addplot [
color=black,
solid,
mark=square,
mark options={solid}
]
table[row sep=crcr]{
44 8.08596416136617e-15\\
82 5.0592338969142e-15\\
164 1.47647841692069e-14\\
320 1.99924447753013e-14\\
642 1.71577822553291e-14\\
1276 1.45433485010627e-14\\
2542 1.93048125477127e-14\\
5080 6.69026253712372e-14\\
};
\addlegendentry{Polygonal SBFEM};

\addplot [
color=black,
solid,
mark=o,
mark options={solid}
]
table[row sep=crcr]{
44 0.025091425583601\\
82 0.013129562762179\\
164 0.007792580402508\\
320 0.006903643309866\\
642 0.005492911102597\\
1276 0.002943804283943\\
2542 0.002550389473212\\
5080 0.002047256824174\\
};
\addlegendentry{Polygonal nSFEM};

\end{axis}
\end{tikzpicture}%
\caption{Equilibrium patch test: Convergence results for the relative error in the displacement norm $(L^2)$.}
\label{fig:plateHoleConveResults}
\end{figure}

\fref{fig:plateHoleConveResults} shows the relative error in displacement norm $L^2$ as a function of total number of degrees of freedom (dofs) for the different formulations. The same mesh is used when comparing the performances of the different polygonal elements. It is seen that the $L^2-$error is very small, ~1e${-8}$ in case of polygonal FEM, but does not decrease down to machine precision. This can be attributed to the integration error as discussed in~\cite{talischipaulino2013}. The order of accuracy is similar to that reported in~\cite{sukumartabarraei2004} for polygonal FEM with Laplace interpolants. In case of the polygonal nSFEM, the method fails to pass the patch test. Again, this can be attributed to the numerical integration, in which the simple averaging technique introduces additional errors when computing the stiffness matrix.  The scaled boundary polygon formulation satisfies the patch test down to machine precision. This can be attributed to the semi-analytical nature of the SBFEM technique. 

Higher order interpolants can be conveniently formulated using the scaled boundary polygon formulation. This involves only increasing the order of the one-dimensional shape function $\mathbf{N}(\eta)$ in \Eref{eqn:dispapprox}. The process is simpler compared with the other polygon formulations reported in the literature. For example, the approach developed by Sukumar~\cite{sukumar2013} requires the solution of an optimization problem. Table \ref{table:simplepatchresult} presents the relative error in the displacement norm $L^2$ for different orders of shape functions along the each edge of the polygonal element. It is seen that the polygon SBFEM satisfies the linear patch test down to machine precision. 

\begin{table}[htpb]
\centering
\caption{Relative error in the displacement norm for the equilibrium patch test.}
\begin{tabular}{lrrr}
\hline
Number of & \multicolumn{3}{c}{Order of shape functions.} \\
\cline{2-4}
Polygons & $p=$ 1 & $p=$ 2 & $p=$ 3 \\
\hline
1 & 1.1955e${-15}$ & 2.4484e${-14}$ & 2.4442e${-14}$\\
10 & 1.1228e${-14}$ & 7.1956e${-14}$ & 3.4006e${-14}$\\
20 & 1.5151e${-14}$ & 5.5886e${-14}$ & 3.9309e${-14}$\\
30 & 1.7593e${-14}$ & 4.7998e${-14}$ & 3.2566e${-14}$\\
\hline
\end{tabular}
\label{table:simplepatchresult}
\end{table}

\subsection{Cantilever beam}
\begin{figure}[htpb]
\centering
\scalebox{0.7}{\input{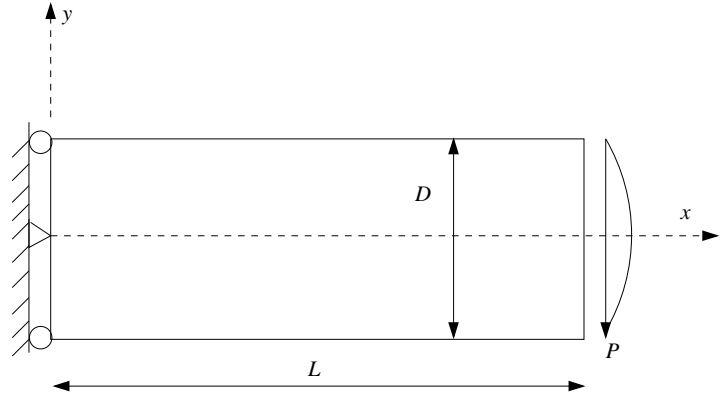}}
\caption{Cantilever beam: Geometry and boundary conditions.}
\label{fig:cantileverfig}
\end{figure}

\begin{figure}[htpb]
\centering
\includegraphics[scale=1]{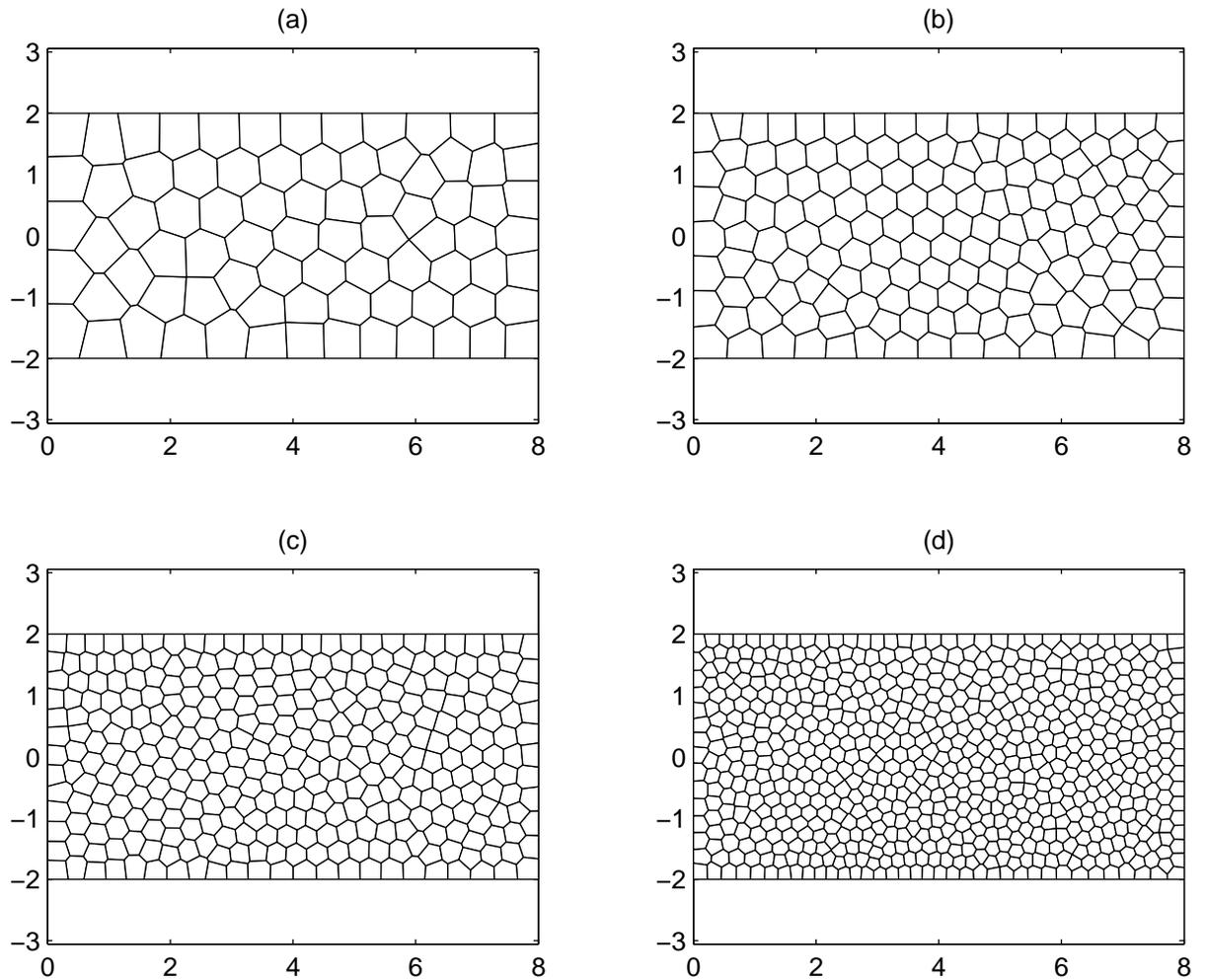}
\caption{Thick cantilever beam: domain discretized with polygonal elements: (a) 80 elements; (b) 160 elements; (c) 320 elements and (d) 640 elements.}
\label{fig:cantmesh}
\end{figure}

A two-dimensional cantilever beam subjected to a parabolic shear load at the free end is examined as shown in \fref{fig:cantileverfig}. The geometry is: length $L=$ 10, height $D=$ 2. The material properties are: Young's modulus, $E=$ 3e$^7$, Poisson's ratio $\nu=$ 0.25 and the parabolic shear force $P=$ 150. The exact solution for displacements are given by:

\begin{align}
u(x,y) &= \frac{P y}{6 \overline{E}I} \left[ (9L-3x)x + (2+\overline{\nu}) \left( y^2 - \frac{D^2}{4} \right) \right] \nonumber \\
v(x,y) &= -\frac{P}{6 \overline{E}I} \left[ 3\overline{\nu}y^2(L-x) + (4+5\overline{\nu}) \frac{D^2x}{4} + (3L-x)x^2 \right]
\label{eqn:cantisolution}
\end{align}
where $I = D^3/12$ is the moment of inertia, $\overline{E} = E$, $\overline{\nu} = \nu$ and $\overline{E} = E/(1-\nu^2)$, $\overline{\nu} = \nu/(1-\nu)$ for plane stress and plane strain, respectively. \fref{fig:cantmesh} shows a sample polygonal mesh used for this study.  

\begin{figure}
\centering
\setlength\figureheight{8cm} 
\setlength\figurewidth{10cm}
\subfigure[]{
%
%
%
%
\begin{tikzpicture}

\begin{axis}[%
width=\figurewidth,
height=\figureheight,
scale only axis,
xmode=log,
xmin=1,
xmax=100,
xminorticks=true,
xlabel={$\text{log((total dof)}^{\text{1/2}}\text{)}$},
ymode=log,
ymin=1e-07,
ymax=1,
yminorticks=true,
ylabel={log(Relative error in displacement)},
legend style={draw=black,fill=white,legend cell align=left},
legend pos={south west}
]
\addplot [
color=black,
solid,
mark=diamond,
mark options={solid}
]
table[row sep=crcr]{
6.6332495807108 0.077930123182295\\
9.05538513813742 0.038128851629516\\
12.7279220613579 0.021569714681425\\
17.832554500127 0.009410147260429\\
25.3377189186399 0.004250971760058\\
35.749125863439 0.002381577687072\\
50.4975246918104 0.001038748189487\\
};
\addlegendentry{Polygonal FEM (m=2.13)};

\addplot [
color=black,
solid,
mark=o,
mark options={solid}
]
table[row sep=crcr]{
6.6332495807108 0.314139280862195\\
9.05538513813742 0.256288820865718\\
12.7279220613579 0.155450120984294\\
17.832554500127 0.140684208020197\\
25.3377189186399 0.086411133305727\\
35.749125863439 0.085120866316996\\
50.4975246918104 0.072439658207499\\
};
\addlegendentry{Polygonal nSFEM (m=0.72)};

\addplot [
color=black,
solid,
mark=square,
mark options={solid}
]
table[row sep=crcr]{
6.6332495807108 0.029801861377447\\
9.05538513813742 0.010841883468604\\
12.7279220613579 0.0077256209097\\
17.832554500127 0.001416139069333\\
25.3377189186399 0.0006452521987667\\
35.749125863439 0.0004476055800628\\
50.4975246918104 0.0001367009735622\\
};
\addlegendentry{Polygonal SBFEM (p=1,m=2.67)};

\addplot [
color=black,
dashed,
mark=square,
mark options={solid}
]
table[row sep=crcr]{
10.295630140987 0.000913378553323\\
14.2126704035519 0.0003911586760146\\
20.0499376557634 0.0001290314067945\\
28.3196045170126 5.50929713292e-05\\
39.9749921826134 1.01363695647e-05\\
56.5508620623948 1.795980228e-06\\
79.9624912068152 5.245280663e-07\\
};
\addlegendentry{Polygonal SBFEM (p=2,m=3.63)};

\end{axis}
\end{tikzpicture}%
}
\subfigure[]{
%
%
%
%
\begin{tikzpicture}

\begin{axis}[%
width=\figurewidth,
height=\figureheight,
scale only axis,
xmode=log,
xmin=1,
xmax=100,
xminorticks=true,
xlabel={$\text{log((total dof)}^{\text{1/2}}\text{)}$},
ymode=log,
ymin=0.0001,
ymax=1,
yminorticks=true,
ylabel={log(Relative error in energy)},
legend style={draw=black,fill=white,legend cell align=left},
legend pos={south west}
]
\addplot [
color=black,
solid,
mark=diamond,
mark options={solid}
]
table[row sep=crcr]{
6.6332495807108 0.324401722151178\\
9.05538513813742 0.227510865808106\\
12.7279220613579 0.164468288504801\\
17.832554500127 0.124816703370485\\
25.3377189186399 0.0866610274595894\\
35.749125863439 0.0580977010029557\\
50.4975246918104 0.0423873435367917\\
};
\addlegendentry{Polygonal FEM (m=1.00)};

\addplot [
color=black,
solid,
mark=o,
mark options={solid}
]
table[row sep=crcr]{
6.6332495807108 0.576204856992045\\
9.05538513813742 0.520930275522621\\
12.7279220613579 0.401193789257011\\
17.832554500127 0.387293232566969\\
25.3377189186399 0.302191497422609\\
35.749125863439 0.29860587936965\\
50.4975246918104 0.273139911079083\\
};
\addlegendentry{Polygonal nSFEM (m=0.4)};

\addplot [
color=black,
solid,
mark=square,
mark options={solid}
]
table[row sep=crcr]{
6.6332495807108 0.242375510342259\\
9.05538513813742 0.155478685572924\\
12.7279220613579 0.113571300992719\\
17.832554500127 0.0841825237204677\\
25.3377189186399 0.0538213157639837\\
35.749125863439 0.0368707975702187\\
50.4975246918104 0.0265768034273381\\
};
\addlegendentry{Polygonal SBFEM (p=1,m=1.09)};

\addplot [
color=black,
dashed,
mark=square,
mark options={solid}
]
table[row sep=crcr]{
6.6332495807108 0.0129149259485574\\
9.05538513813742 0.0074721852995756\\
12.7279220613579 0.00467455469535734\\
17.832554500127 0.00365231837225664\\
25.3377189186399 0.00194615334227347\\
35.749125863439 0.00091094739442808\\
50.4975246918104 0.000458056719406388\\
};
\addlegendentry{Polygonal SBFEM (p=2,m=1.63)};

\end{axis}
\end{tikzpicture}%
}
\caption{Bending of thick cantilever beam: Convergence results for (a) the relative error in the displacement norm $(L^2)$ and (b) the relative error in the energy norm. The rate of convergence is also shown, where $m$ is the average slope. In case of the polygonal SBFEM, $p$ denotes the order of the shape functions along each edge of the polygon.}
\label{fig:cantiConveResults}
\end{figure}
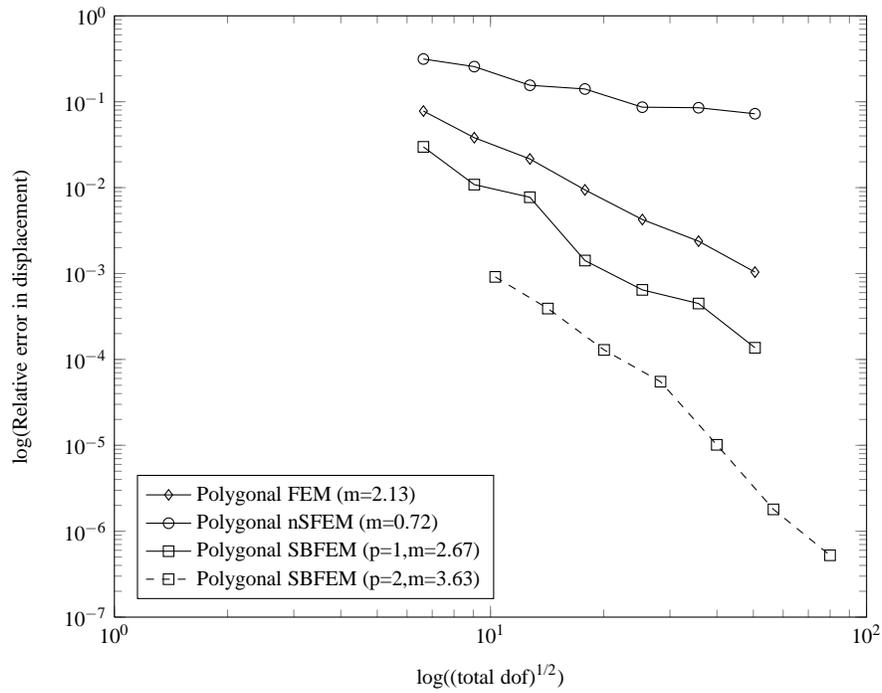
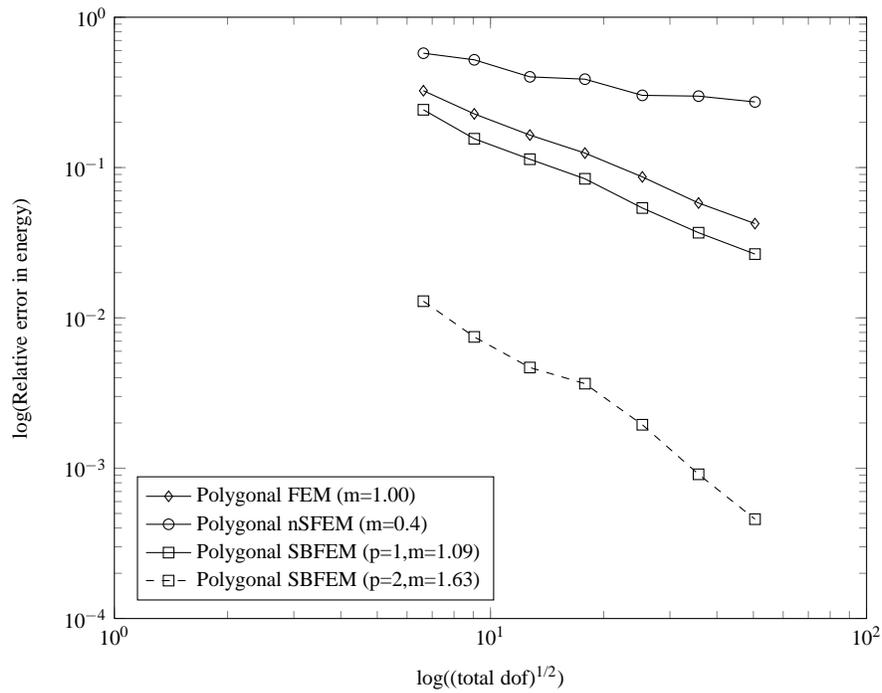

The numerical convergence of the relative error in the displacement norm and the relative error in the energy norm is shown in \fref{fig:cantiConveResults}. The results from different approaches are compared with the available analytical solution. The numerical integration error is present in both the Polygonal FEM and the Polygonal nSFEM, but is more pronounced in case of the Polygonal nSFEM. Both the Polygonal FEM and the Polygonal SBFEM yields optimal convergence in $L^2$ and $H^1$, whilst the polygonal nSFEM yields sub-optimal convergence rate. This can be attributed to the error in the numerical integration introduced by the simple averaging technique. It is seen that with mesh refinement, all the methods converge to the exact solution. An estimation of the convergence rate is also shown. From \fref{fig:cantiConveResults}, it can be observed that the scaled boundary polygon formulation yields more accurate results. It is also noted that in case of the polygonal SBFEM, no special numerical integration technique is required to compute the stiffness matrix and extension to higher order shape function is straight forward. \fref{fig:cantiConveResults} also shows the error in the displacement norm when quadratic shape functions are used along each edge within the scaled boundary polygon formulation. It is seen that as the order of the shape functions is increased, the error decreases while the convergence rate increases.

\subsection{Infinite plate with a circular hole}
\begin{figure}[htpb]
\centering
\includegraphics[scale=1]{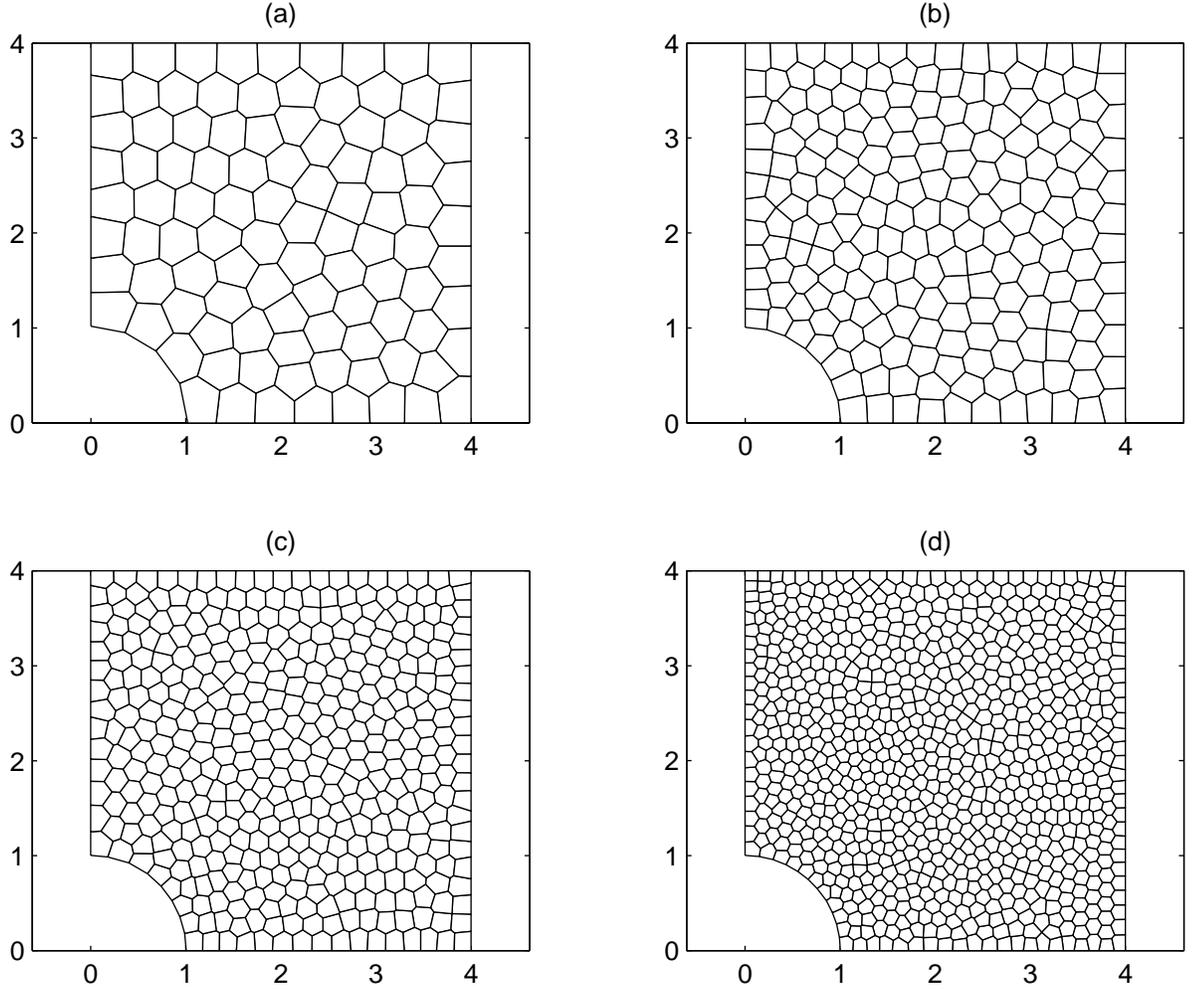}
\caption{Plate with a circular hole: domain discretized with polygonal elements: (a) 100 elements; (b) 200 elements; (c) 400 elements and (d) 800 elements.}
\label{fig:phwmesh}
\end{figure}

In this example, consider an infinite plate with a traction free hole under uniaxial tension $(\sigma=$1$)$ along $x-$axis~\fref{fig:twocrkhole}. The exact solution of the principal stresses in polar coordinates $(r,\theta)$ is given by:
\begin{eqnarray}
\sigma_{11}(r,\theta) &= 1 - \frac{a^2}{r^2} \left( \frac{3}{2} (\cos 2\theta + \cos 4\theta) \right) + \frac{3a^4}{2r^4} \cos 4\theta \nonumber \\
\sigma_{22}(r,\theta) &= -\frac{a^2}{r^2} \left( \frac{1}{2}(\cos 2\theta - \cos 4\theta) \right) - \frac{3a^4}{2r^4} \cos 4\theta \nonumber \\
\sigma_{12}(r,\theta) &= -\frac{a^2}{r^2} \left( \frac{1}{2}(\sin 2\theta + \sin 4\theta) \right) + \frac{3a^4}{2r^4} \sin 4\theta
\end{eqnarray}
where $a$ is the radius of the hole. Owing to symmetry, only one quarter of the plate is modeled. \fref{fig:phwmesh} shows a typical polygonal mesh used for the study. The material properties are: Young's modulus $E=$ 10$^5$ and Poisson's ratio $\nu=$ 0.3. In this example, analytical tractions are applied on the boundary. The domain is discretized with polygonal elements and along each edge of the polygon, the shape function is linear. The convergence rate in terms of the displacement norm is shown in \fref{fig:plateHoleConveResults}. It can be seen that the Polygonal SBFEM yields more accurate results when compared with the Polygonal FEM and the Polygonal nSFEM. The convergence rate of the polygonal FEM and the polygonal SBFEM are similar. The convergence rate of the polygonal nSFEM is inferior compared with the polygonal FEM and the polygonal SBFEM. All the techniques converge to exact energy with mesh refinement. The relative error in terms of the displacement norm for the polygon SBFEM with quadratic shape functions is shown in \fref{fig:plateHoleConveResults}. It is observed that the error decreases as the order of the shape functions is increased.

\begin{figure}[htpb]
\centering
\scalebox{0.9}{\input{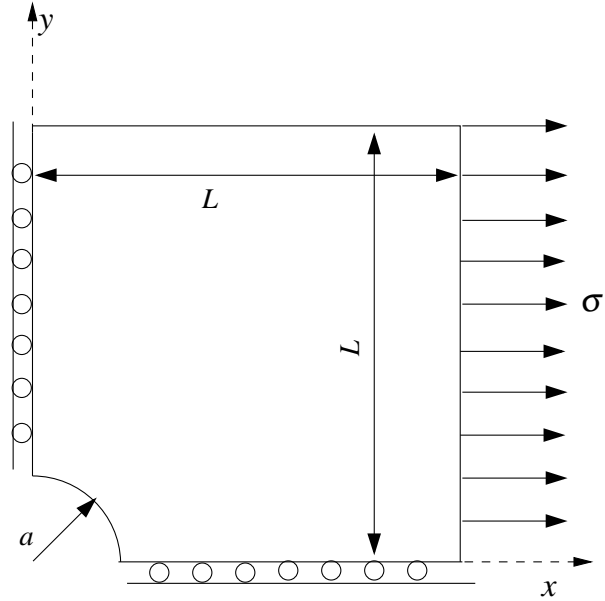}}
\caption{Infinite plate with a circular hole.}
\label{fig:twocrkhole}
\end{figure}

\begin{figure}
\centering
\setlength\figureheight{8cm} 
\setlength\figurewidth{10cm}
%
%
%
%
\begin{tikzpicture}

\begin{axis}[%
width=\figurewidth,
height=\figureheight,
scale only axis,
xmode=log,
xmin=100,
xmax=10000,
xminorticks=true,
xlabel={$\text{log((total dof)}^{\text{1/2}}\text{)}$},
ymode=log,
ymin=1e-06,
ymax=0.1,
yminorticks=true,
ylabel={log(Relative error in displacement)},
legend style={draw=black,fill=white,legend cell align=left},
legend pos={south west}
]
\addplot [
color=black,
solid,
mark=diamond,
mark options={solid}
]
table[row sep=crcr]{
404 0.007201360693307\\
790 0.004368235117637\\
1598 0.00186467533149\\
3178 0.0008014627252884\\
6360 0.000455433128974\\
};
\addlegendentry{Polygonal FEM (m=1.99)};

\addplot [
color=black,
solid,
mark=o,
mark options={solid}
]
table[row sep=crcr]{
404 0.021742978488148\\
790 0.012356149377212\\
1598 0.009242667275117\\
3178 0.00566477327497\\
6360 0.004641240777099\\
};
\addlegendentry{Polygonal nSFEM (m=1.13)};

\addplot [
color=black,
solid,
mark=square,
mark options={solid}
]
table[row sep=crcr]{
404 0.002870985039609\\
790 0.002245866916675\\
1598 0.000658899213228\\
3178 0.0002492743511912\\
6360 0.000142721361221\\
};
\addlegendentry{Polygonal SBFEM (p= 1,m=2.16)};

\addplot [
color=black,
dashed,
mark=square,
mark options={solid}
]
table[row sep=crcr]{
1002 0.000301\\
2002 0.000113\\
3982 2.16e-05\\
7970 8.35e-06\\
};
\addlegendentry{Polygonal SBFEM (p=2,m=3.46)};

\end{axis}
\end{tikzpicture}%
\caption{Infinite plate with a circular hole: Convergence results for the relative error in the displacement norm $(L^2)$. The rate of convergence is also shown, where $m$ is the average slope.}
\label{fig:plateHoleConveResults}
\end{figure}
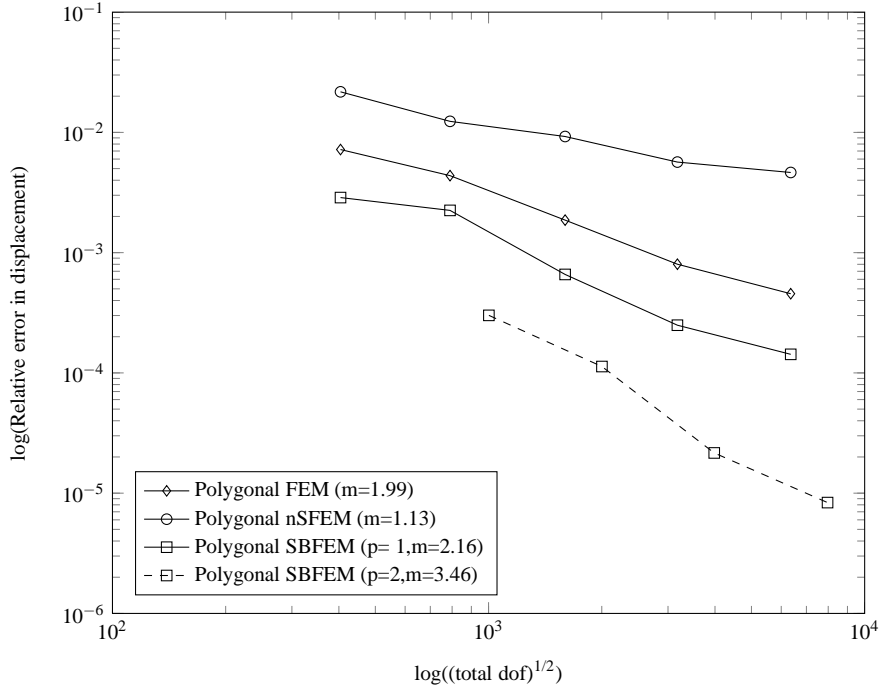

\subsection{Application to linear elastic fracture mechanics}
The polygonal finite element techniques discussed above can be applied to problems in linear elastic fracture mechanics. However, to accurately capture the singularity at the crack tip, a very fine mesh in combination with singular elements at the crack tip is usually required. This poses numerical difficulties when the crack evolves. Another possibility is to augment the existing finite element approximation basis with additional functions~\cite{tabarraeisukumar2008}. In the literature, the latter method is referred to as the `extended finite element method' (XFEM). The flexibility provided by the XFEM, also has certain difficulties,  such as numerical integration and blending of enriched and non-enriched elements~\cite{friesbelytschko2010}. In this section, we discuss the application of the polygonal SBFEM to problems with strong discontinuities and singularities. By exploiting the special characteristics of the scaling center~\cite{wolfsong2001}, the stress intensity factors can be computed directly.  When modeling a cracked/notched structure, the scaling centre is placed at the crack tip. The straight crack/notch edges are formed by scaling the nodes $A$ and $B$ on the boundary and are not discretized (see ~\fref{fig:Polygon-representation-by}). In this study, a polygonal mesh that conforms to the crack surface is generated.

\paragraph{Calculation of the stress intensity factors} A unique feature of the scaled boundary polygon formulation is that stress singularities of any kind e.g. in cracks notches and material junctions, if present, are analytically represented in the radial displacement functions $\mathbf{u}(\xi)$. When a crack is modelled by a polygon element with the scaling centre chosen at its crack tip in~\fref{fig:Polygon-representation-by}, some of the eigenvalues in $\boldsymbol{\Lambda}_{\mathrm{n}}$ satisfy $-1<\lambda\left(\boldsymbol{\Lambda}_{\mathrm{n}}\right)<0$. The stress field $\boldsymbol{\sigma}(\xi,\eta)$ can be expressed as~\cite{wolfsong2001}:

\begin{align}
\boldsymbol{\sigma}(\xi,\eta)= & \boldsymbol{\Psi}_{\sigma}(\eta)\xi^{-\boldsymbol{\Lambda}_{\mathrm{n}}-\mathbf{I}}\mathbf{c}
\label{eq:stress field complete}
\end{align}
where $\boldsymbol{\Lambda}_{\mathrm{n}}=\mathrm{diag}\left(\lambda_{1},\,\lambda_{2},\,...,\lambda_{n}\right)$ are the eigenvalues, $\mathbf{c} = \boldsymbol{\Phi}^{-1}_{\mathrm{u}} \uu_b$, $\uu_b$ is the displacement vector on the boundary $\xi=$ 1, obtained by solving the algebraic system of equations. The stress mode $\boldsymbol{\Psi}_{\sigma}(\eta)$ is defined as
\begin{align}
\boldsymbol{\Psi}_{\sigma}(\eta)= & \mathbf{D}\left(-\mathbf{B}_{1}(\eta)\boldsymbol{\Phi}_{\mathrm{u}}\boldsymbol{\Lambda}_{\mathrm{n}}+\mathbf{B}_{2}(\eta)\boldsymbol{\Phi}_{\mathrm{u}}\right)\label{eq:stress mode}
\end{align}
It can be discerned from~\Eref{eq:stress field complete} that these eigenvalues lead to singular stresses at the crack tip. This enables the stress intensity factors (SIF) to be computed directly from their definitions. For a crack that is aligned with the Cartesian coordinate axes shown in \fref{fig:crkPolySBFEM}, the SIFs are defined as:
\begin{figure}[htpb]
\centering
\scalebox{0.7}{\input{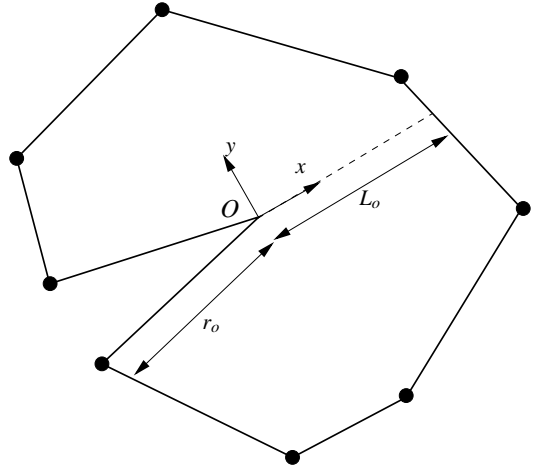}}
\caption{A cracked domain modelled by SBFEM and the definition of local coordinate system, where the `black' dots represent the nodes.}
\label{fig:crkPolySBFEM}
\end{figure}

\begin{equation}
\left\{ \begin{array}{c}
K_{\mathrm{I}}\\
K_{\mathrm{II}}
\end{array}\right\} =\lim_{r\rightarrow0}\left\{ \begin{array}{c}
\sqrt{2\pi r}\left.\sigma_{yy}\right|_{\theta=0}\\
\sqrt{2\pi r}\left.\tau_{xy}\right|_{\theta=0}
\end{array}\right\} \label{eq:SIF-stress}
\end{equation}
In the limit when $r\rightarrow0$ or $\xi\rightarrow0$,~\Eref{eq:stress field complete} shows that the stresses vanish when the real parts of eigenvalues in $\boldsymbol{\Lambda}_{\mathrm{n}}$, $\lambda\left(\boldsymbol{\Lambda}_{\mathrm{n}}\right)\leq-1$. The diagonal submatrix $\boldsymbol{\Lambda}_{\mathrm{n}}^{(\mathrm{s})}=-0.5\mathbf{I}$ becomes dominant. For this diagonal submatrix, the singular stress field $\boldsymbol{\sigma}^{(\mathrm{s})}(\xi,\eta)$ can be expressed as: 
\begin{equation}
\begin{array}{ccc}
\boldsymbol{\sigma}^{(\mathrm{s})}(\xi,\eta) & = & \boldsymbol{\Psi}_{\sigma}^{(\mathrm{s})}\left(\eta\left(\theta\right)\right)\xi^{-0.5\mathbf{I}}\mathbf{c}^{(\mathrm{s})}\end{array}\label{eq:Singular stress field}
\end{equation}
 where $\boldsymbol{\Psi}_{\sigma}^{(\mathrm{s})}\left(\eta\left(\theta\right)\right)=\left[\begin{array}{ccc}
\boldsymbol{\Psi}_{\sigma_{xx}}^{(\mathrm{s})}\left(\eta\left(\theta\right)\right) & \boldsymbol{\Psi}_{\sigma_{yy}}^{(\mathrm{s})}\left(\eta\left(\theta\right)\right) & \boldsymbol{\Psi}_{\tau_{xy}}^{(\mathrm{s})}\left(\eta\left(\theta\right)\right)\end{array}\right]^{\mathrm{T}}$ and $\mathbf{c}^{(\mathrm{s})}$ are the stress modes and integration constants corresponding to the diagonal submatrix in $\boldsymbol{\Lambda}_{\mathrm{n}}$ satisfying $-1<\lambda\left(\boldsymbol{\Lambda}_{\mathrm{n}}\right)<0$. Substituting the stress components in~\Eref{eq:Singular stress field} at angle $\theta=0$ into~\Eref{eq:SIF-stress} and using the relation $\xi=\frac{r}{L_{0}}$ at $\theta=$0, the SIF is
\begin{align}
\left\{ \begin{array}{c}
K_{\mathrm{I}}\\
K_{\mathrm{II}}
\end{array}\right\} = & \sqrt{2\pi L_{0}}\left\{ \begin{array}{c}
\boldsymbol{\Psi}_{\sigma_{yy}}^{(s)}(\eta(\theta=0))\mathbf{c}^{(s)}\\
\boldsymbol{\Psi}_{\tau_{xy}}^{(s)}(\eta(\theta=0))\mathbf{c}^{(s)}
\end{array}\right\} 
\end{align}

\subsubsection{Plate with double edge crack in tension}
\begin{figure}[htpb]
\centering
\subfigure[]{\scalebox{0.7}{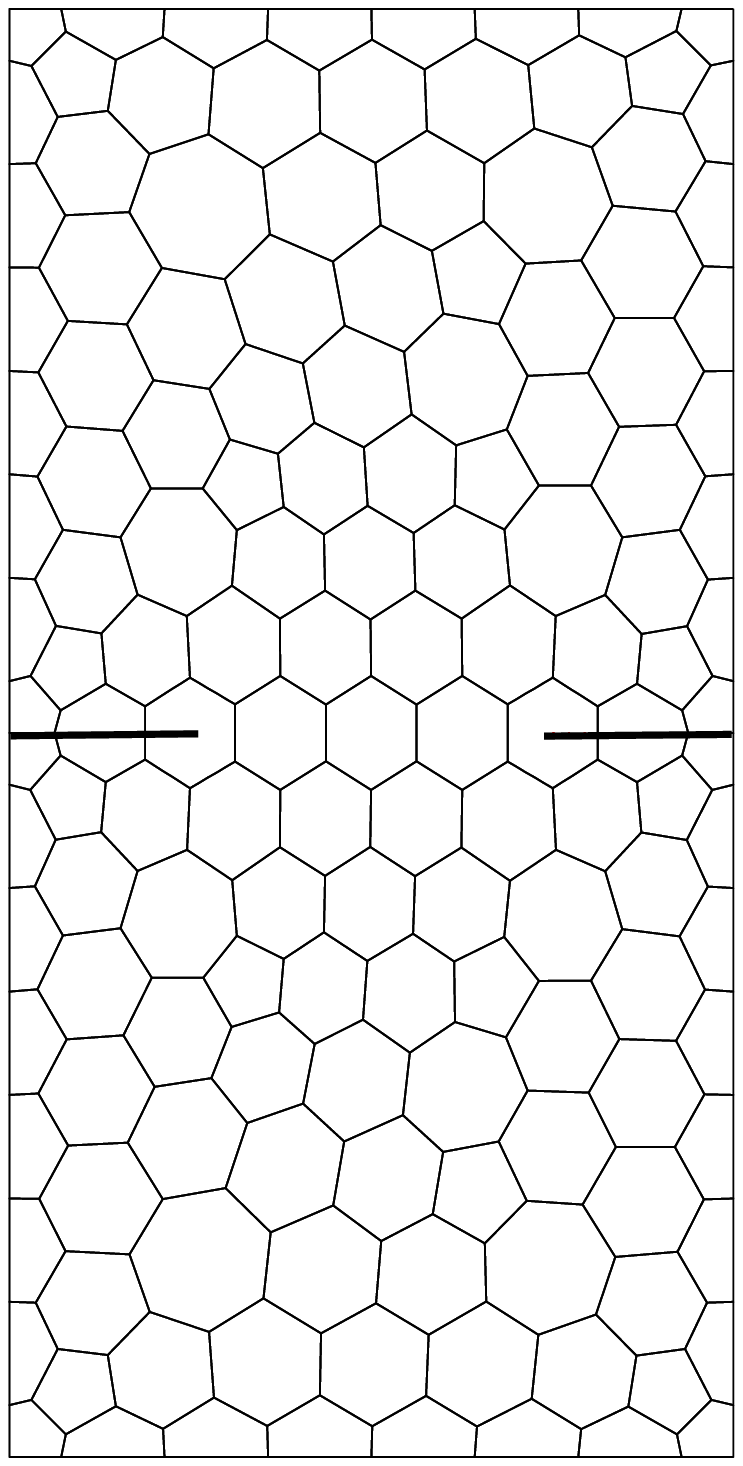}}
\subfigure[]{\includegraphics[scale=0.45]{./Figures/dedgcrk}}
\caption{Plate with an edge under tension: (a) geometry and boundary conditions and (b) domain discretized with polygonal elements.}
\label{fig:edgcrkgeom}
\end{figure}

The plate with double edge crack subjected to a uniform tension at both ends as shown in \fref{fig:edgcrkgeom} is considered. In the computations, the ratio of the crack length, $a$, to the width of the plate, $H$, is $a/H=$ 0.25. The material properties of the plate are: Young's modulus, $E=$ 200 GPa and Poisson's ratio, $\nu=0.3$. In this example, plane stress conditions are assumed. The empirical mode I SIF that is given by:
\begin{equation}
K_{\rm{I}}^{\mathrm{ref}}=C\sigma\sqrt{\pi a}\label{eq: KI ref}
\end{equation}
where $C$ is a correction factor. For $a/b > 0.4, b = H/2$, the correction factor is given by~\cite{tadaparis2000}:
\begin{equation}
C=1.12+0.203\left(\frac{a}{b}\right)-1.197\left(\frac{a}{b}\right)^{2}+1.930\left(\frac{a}{b}\right)^{3} 
\label{eqn:correctionfactor1}
\end{equation}

The above correction factor is for an infinite plate and with an accuracy of 2\%~\cite{tadaparis2000}. For the chosen parameters, the reference normalized SIF is, $K_{\rm I}/\sqrt{\pi a} = $ 1.1635. The plate discretized using a polygon mesh of arbitrary number of sides. For the cracked polygon, each edge is discretized with twelve elements so that the angular variation of the SIF can be captured accurately. Table \ref{table:dedgcrkmeshcon} presents the results from the Polygonal SBFEM with linear and quadratic shape functions along the each edge of the element. The results from the Polygonal SBFEM are compared with the XFEM results~\cite{tabarraeisukumar2008}. The number of nodes is maintained approximately the same as reported in \cite{tabarraeisukumar2008}. It can be seen that the results from the Polygonal SBFEM are in good agreement with the results in the literature. It is noted that, no special functions are required in case of the polygonal SBFEM and the SIF can be computed directly as outlined above. In Table \ref{table:dedgcrkmeshconvehl3}, the results from the Polygonal SBFEM are compared with the Body Force Method (BFM)~\cite{saimotomotomura2010} with $L/H=$ 3. It is observed that the results from the Polygonal SBFEM are in very good agreement with the BFM and also with the empirical relation given in~\cite{tadaparis2000}. Also, the error decreases as the order of the shape functions is increased.

\begin{table}[htpb]
\centering
\caption{Normalized mode I SIF for double-edge cracked plate in tension for $L/H=$ 2, $a/H=$ 0.25.}
\begin{tabular}{lrrrrrrrr}
\hline 
&\multicolumn{2}{c}{Ref.~\cite{tabarraeisukumar2008}} && \multicolumn{2}{c}{Polygonal SBFEM ($p=$ 1)} && \multicolumn{2}{c}{Polygonal SBFEM ($p=$ 2)} \\
\cline{2-3} \cline{5-6} \cline{8-9}
& Number of nodes & $K_{\rm I}/\sqrt{\pi a}$ && Number of nodes & $K_{\rm I}/\sqrt{\pi a}$ && Number of nodes & $K_{\rm I}/\sqrt{\pi a}$ \\
\hline
Mesh 1& - & - && 116 & 1.1582 && 255 & 1.1702 \\
Mesh 2 & - & - && 148 & 1.1599 && 381 & 1.1702 \\
Mesh 3& 402 & 1.1622 && 400 & 1.1662 && 946 & 1.1703 \\
Mesh 4&1002 & 1.1625 && 1050 & 1.1679 && 2549 & 1.1703 \\
Mesh 5&2001 & 1.1631 && 2012 & 1.1679 && 4934 & 1.1703 \\
\hline 
\end{tabular}
\label{table:dedgcrkmeshcon}
\end{table}

\begin{table}[htpb]
\centering
\caption{Normalized mode I SIF for double-edge cracked plate in tension for $L/H=$ 3, $a/H=$ 0.25.}
\begin{tabular}{lrrr}
\hline 
Number of Polygons && \multicolumn{2}{c}{Polygonal SBFEM, $K_{\rm I}/\sqrt{\pi a}$} \\
\cline{3-4}
&& $p=$1 & $p=$ 2 \\
\hline
194 && 1.16531 &  1.16925 \\
637 && 1.16680 &  1.16927 \\
1327 && 1.16843 & 1.16927 \\
Ref.~\cite{saimotomotomura2010} && 1.16925 & 1.16925 \\
\hline 
\end{tabular}
\label{table:dedgcrkmeshconvehl3}
\end{table}

\subsubsection{Inclined crack in plate}
\begin{figure}[htpb]
\centering
\scalebox{0.7}{\input{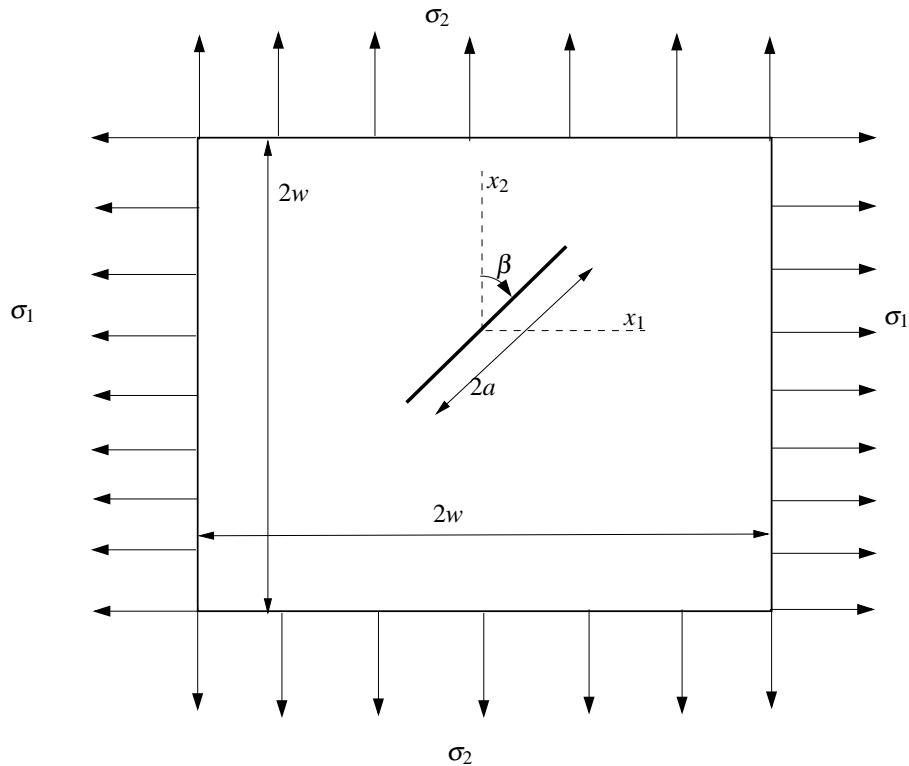}}
\caption{Plate with an oblique crack: geometry and boundary conditions.}
\label{fig:inclcrkgeom}
\end{figure}

In the next example, a plate with an angled crack subjected to far field bi-axial stress field, $\bvsig$ (see \fref{fig:inclcrkgeom} with $a/w=$ 0.1, $\sigma_1=$ 1 and $\sigma_2=$ 2 is considered. In this example, the mode I and the mode II SIFs, $K_{\rm I}$ and $K_{\rm II}$, respectively, are obtained as a function of the crack angle $\beta$. For the loads shown, the analytical SIF for an infinite plate are given by~\cite{aliabadirooke1987,tabarraeisukumar2008}:
\begin{align}
K_{\rm I} &= (\sigma_2 \sin^2 \beta + \sigma_1 \cos^2 \beta) \sqrt{\pi a} \nonumber \\
K_{\rm II} &= (\sigma_2 - \sigma_1) \sin \beta \cos \beta \sqrt{\pi a}
\label{eqn:inclcrcempirical}
\end{align}

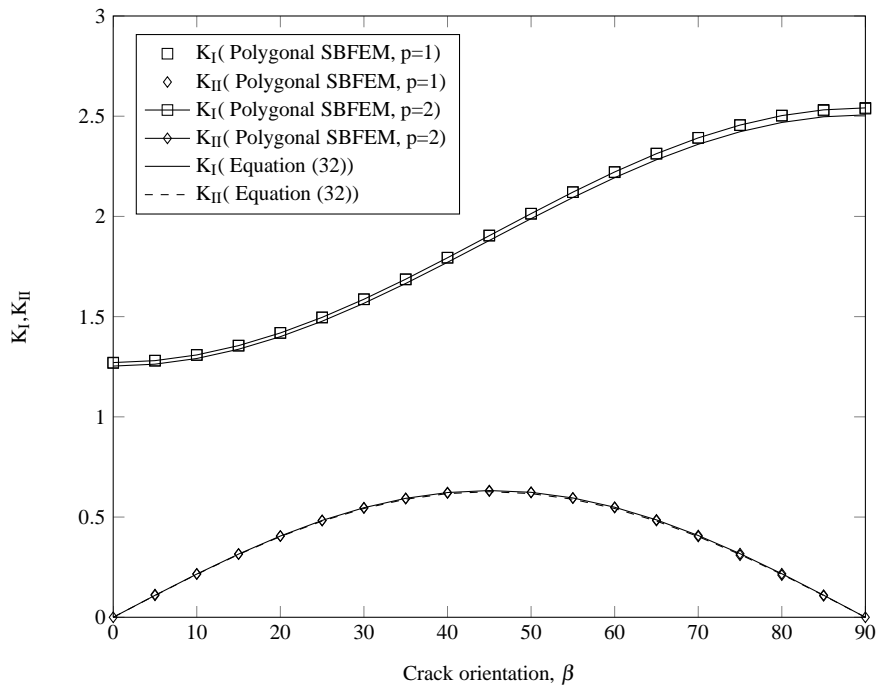
\begin{figure}
\centering
\setlength\figureheight{8cm} 
\setlength\figurewidth{10cm}
%
%
%
%
\begin{tikzpicture}

\begin{axis}[%
width=\figurewidth,
height=\figureheight,
scale only axis,
xmin=0,
xmax=90,
xlabel={$\text{Crack orientation, }\beta$},
ymin=0,
ymax=3,
ylabel={$\text{K}_\text{I}\text{,K}_{\text{II}}$},
legend style={draw=black,fill=white,legend cell align=left},
legend pos = north west
]
\addplot [
color=black,
only marks,
mark=square,
mark options={solid}
]
table[row sep=crcr]{
90 2.537759799\\
85 2.528142379\\
80 2.502932918\\
75 2.45530286\\
70 2.392371696\\
65 2.313620309\\
60 2.220332028\\
55 2.120963658\\
50 2.013276562\\
45 1.905161514\\
40 1.794475375\\
35 1.685161535\\
30 1.586560236\\
25 1.495524323\\
20 1.419225557\\
15 1.35525562\\
10 1.30884275\\
5 1.280201984\\
0 1.268877054\\
};
\addlegendentry{$\text{K}_\text{I}\text{( Polygonal SBFEM, p=1)}$};

\addplot [
color=black,
only marks,
mark=diamond,
mark options={solid}
]
table[row sep=crcr]{
90 -1.9084e-05\\
85 0.107046974\\
80 0.210216247\\
75 0.30959648\\
70 0.40302619\\
65 0.482487366\\
60 0.547321699\\
55 0.594739703\\
50 0.623085351\\
45 0.629512855\\
40 0.619817332\\
35 0.591892628\\
30 0.544954307\\
25 0.482619373\\
20 0.405294887\\
15 0.315982054\\
10 0.216764833\\
5 0.113373057\\
0 -9.541e-06\\
};
\addlegendentry{$\text{K}_{\text{II}}\text{( Polygonal SBFEM, p=1)}$};

\addplot [
color=black,
solid,
mark=square,
mark options={solid}
]
table[row sep=crcr]{
90 2.541523015\\
85 2.531724231\\
80 2.503045401\\
75 2.455826475\\
70 2.391894053\\
65 2.313061818\\
60 2.221833157\\
55 2.121101602\\
50 2.013833626\\
45 1.903935177\\
40 1.793673323\\
35 1.68684585\\
30 1.586772271\\
25 1.496407491\\
20 1.418510582\\
15 1.355360898\\
10 1.308838633\\
5 1.280419832\\
0 1.27075924\\
};
\addlegendentry{$\text{K}_\text{I}\text{( Polygonal SBFEM, p=2)}$};

\addplot [
color=black,
solid,
mark=diamond,
mark options={solid}
]
table[row sep=crcr]{
90 1.81e-07\\
85 0.110291207\\
80 0.217272441\\
75 0.317498938\\
70 0.408010257\\
65 0.485994539\\
60 0.549102141\\
55 0.595397499\\
50 0.623520915\\
45 0.632347708\\
40 0.622414096\\
35 0.59354834\\
30 0.546720005\\
25 0.483341116\\
20 0.405395292\\
15 0.31521964\\
10 0.215576098\\
5 0.109435875\\
0 9e-08\\
};
\addlegendentry{$\text{K}_{\text{II}}\text{( Polygonal SBFEM, p=2)}$};

\addplot [
color=black,
solid
]
table[row sep=crcr]{
90 2.5066\\
85 2.4971\\
80 2.4688\\
75 2.4226\\
70 2.36\\
65 2.2828\\
60 2.1933\\
55 2.0943\\
50 1.9888\\
45 1.88\\
40 1.7711\\
35 1.6656\\
30 1.5666\\
25 1.4771\\
20 1.4\\
15 1.3373\\
10 1.2911\\
5 1.2628\\
0 1.2533\\
};
\addlegendentry{$\text{K}_\text{I} $( \Eref{eqn:inclcrcempirical})};

\addplot [
color=black,
dashed
]
table[row sep=crcr]{
90 0\\
85 0.1088\\
80 0.2143\\
75 0.3133\\
70 0.4028\\
65 0.48\\
60 0.5427\\
55 0.5888\\
50 0.6171\\
45 0.6266\\
40 0.6171\\
35 0.5889\\
30 0.5427\\
25 0.48\\
20 0.4028\\
15 0.3133\\
10 0.2143\\
5 0.1088\\
0 0\\
};
\addlegendentry{$\text{K}_\text{II} $( \Eref{eqn:inclcrcempirical})};

\end{axis}
\end{tikzpicture}%
\caption{Mode I and mode II SIF for a plate with an inclined crack as a function of crack orientation $\beta$ for $a/w=$ 10.}
\label{fig:inclsifvsangle}
\end{figure}

The material properties of the plate are: Young's modulus, $E=$ 200 GPa and Poisson's ratio, $\nu=0.3$. In this example, the plate is discretized with polygon meshes containing approximately 300 polygons for all values of beta considered. For the element containing the crack tip, five elements (linear or quadratic) were used along each edge to capture the angular variation of the SIF accurately. The variation of the mode I and the mode II SIF with crack orientation $\beta$ are presented in Table \ref{table:inclinedcracksifs} and \fref{fig:inclsifvsangle}. The results from the XFEM~\cite{tabarraeisukumar2008} is also listed in the table, which is for the case of $a/w=$ 10.  The influence of the order of the shape functions and the plate width-to-the crack length ratio is studied using the polygonal SBFEM. It is observed the results compared very well with the reference solutions. Table \ref{table:influenceofaw} presents the influence of the $a/w$ on the numerical SIF in case of the Polygonal SBFEM with quadratic elements. It is seen that, as the plate width-to-the crack length ratio is increased, the numerical SIFs approach the analytical solutions (see \Eref{eqn:inclcrcempirical}) for the infinite plate.

\begin{table}[htpb]
\centering
\renewcommand{\arraystretch}{1.5}
\caption{Mode I and Mode II SIF for a plate with an inclined crack.}
\begin{tabular}{llrrr}
\hline 
SIF & Method & \multicolumn{3}{c}{$\beta$} \\
\cline{3-5}
 & & 90$^\circ$ & 45$^\circ$ &  0$^\circ$\\
\hline 
\multirow{6}{*}{$K_{\rm I}$} & \Eref{eqn:inclcrcempirical} & 2.5066 & 1.8800 & 1.2533 \\
& Ref.~\cite{tabarraeisukumar2008} & 2.5171 & 1.8782 & 1.2585 \\
& Polygonal SBFEM ($p=$ 1) & 2.5378 & 1.9052 & 1.2689 \\
& Polygonal SBFEM ($p=$ 2) & 2.5415 & 1.9039 & 1.2708 \\
& Polygonal SBFEM ($p=$ 3) & 2.5416 & 1.9039 & 1.2708 \\
& Polygonal SBFEM ($p=$ 4) & 2.5416 & 1.9039 & 1.2708 \\
\cline{2-5}
\multirow{6}{*}{$K_{\rm {II}}$} & \Eref{eqn:inclcrcempirical} & 0.0000 & 0.6266 & 0.0000\\
& Ref.~\cite{tabarraeisukumar2008} & 0.0000 & 0.6260 & 0.0000 \\
& Polygonal SBFEM ($p=$ 1) & 0.0000 & 0.6295 & 0.0000 \\
& Polygonal SBFEM ($p=$ 2) & 0.0000 & 0.6323 & 0.0000 \\
& Polygonal SBFEM ($p=$ 3) & 0.0000 & 0.6324 & 0.0000 \\
& Polygonal SBFEM ($p=$ 4) & 0.0000 & 0.6324 & 0.0000 \\
\hline
\end{tabular}
\label{table:inclinedcracksifs}
\end{table}

\begin{table}[htpb]
\centering
\renewcommand{\arraystretch}{1.5}
\caption{Mode I and Mode II SIF for a plate with an inclined crack.}
\begin{tabular}{lllrr}
\hline 
SIF & $\beta$ & \Eref{eqn:inclcrcempirical} & \multicolumn{2}{c}{Polygonal SBFEM ($p=$ 2)} \\
\cline{4-5}
 &  & Infinite Plate &  $a/w=$ 10 & $a/w=$ 50 \\
 \hline
 \multirow{3}{*}{$K_{\rm I}$} & 90$^\circ$ & 2.5066 & 2.5416 & 2.5077 \\
 & 45$^\circ$ & 1.8800 & 1.9039 &  1.8810 \\
 & 0$^\circ$ & 1.2533 & 1.2708 & 1.2539 \\
 \cline{2-5}
 \multirow{3}{*}{$K_{\rm {II}}$} & 90$^\circ$ & 0.0000 & 0.0000 & 0.0000 \\
 & 45$^\circ$ & 0.6266 & 0.6324 & 0.6265 \\
 & 0$^\circ$ & 0.0000 & 0.0000 & 0.0000\\
\hline
\end{tabular}
\label{table:influenceofaw}
\end{table}

\subsubsection{Two cracks emanating from a hole in a finite plate}
As a last example, aplate with a centrally located hole of radius $r$ and two cracks emanating from a hole (see \fref{fig:twocrkhole}) is considered. The dimensions of the plate are: $H = 2W$ and the radius of the hole $r=2W$. The normalized stress intensity factors are defined as: $F_I = K_I/\sigma \sqrt{\pi a}$ and $F_{II} = K_{II}/\sigma \sqrt{\pi a}$. The results from the polygonal SBFEM are compared with the results of Woo \textit{et al.,}~\cite{woowang1989} and Daux \textit{et al.,}~\cite{dauxmoes2000}. In this example, 1D quadratic element is used along each edge of the polygonal element that contains the crack tip. Table \ref{tab:twocrkholeresu} shows the computed SIFs $F_I$ and $F_{II}$ for several ratios of $a/W$ and for different angles $\theta$. It can be seen from Table \ref{tab:twocrkholeresu}, that the results from the polygonal SBFEM are in very good agreement with the results in the literature. It is noted that, when compared to Daux \textit{et al.,}~\cite{dauxmoes2000}, no additional enrichment functions are required to resolve the singularity at the crack tip. Moreover, higher order elements can be employed for the element containing the crack tip. The scaled boundary polygon formulation when combined with the extended FEM, eliminates the need to compute the enrichment functions and improves the accuracy of the computed SIF. Interested readers are referred to~\cite{natarajansong2013} where the XFEM was combined with the SBFEM for quadrilateral elements.

\begin{figure}[htpb]
\centering
\subfigure[]{\scalebox{0.7}{\input{./Figures/holecrack.pstex_t}}}
\subfigure[]{\includegraphics[scale=0.45]{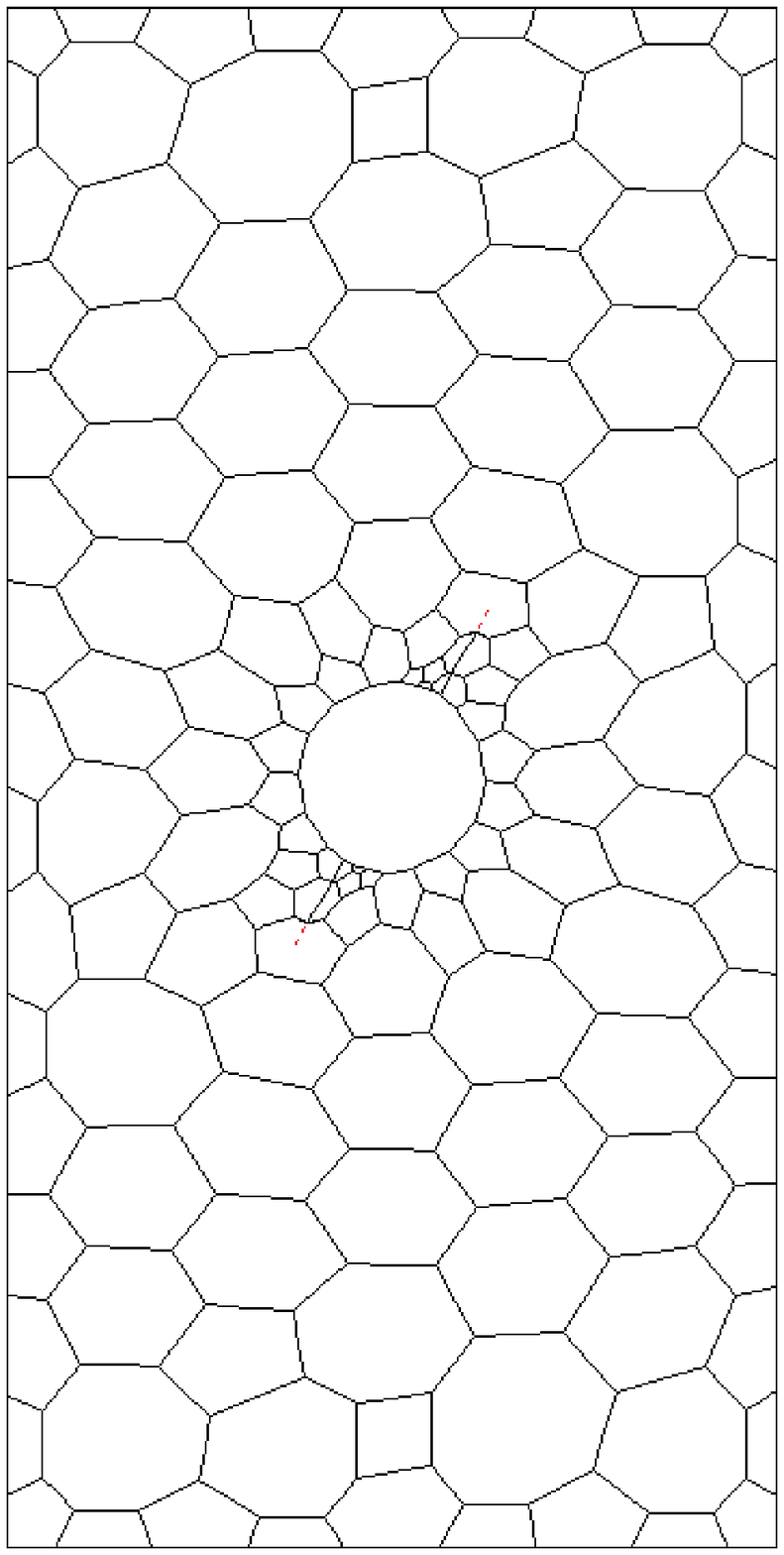}}
\caption{Two cracks emanating from a hole in a finite plate : (a) geometry and boundary conditions and (b) domain discretized with polygonal elements.}
\label{fig:twocrkhole}
\end{figure}

\begin{table}
\centering
\renewcommand{\arraystretch}{1.5}
\caption{Normalised stress intensity factors for two cracks emanating from a hole with various ratios $a/W$ and for different angles $\theta$.}
\begin{tabular}{cllccccccc}
\hline 
 &  &  & \multicolumn{7}{c}{$a/W$}\tabularnewline
\cline{4-10} 
$\theta$ &  &  & 0.3 & 0.4 & 0.5 & 0.6 & 0.7 & 0.8 & 0.9\\
\hline 
\multirow{3}{*}{$0$} &  & Polygonal SBFEM ($p=$ 2) & 1.090 & 1.214 & 1.283 & 1.394 & 1.577 & 1.902 & 2.639\\
 & $F_{I}$ & Ref.~\cite{dauxmoes2000} & 1.082 & 1.207 & 1.286 & 1.389 & 1.558 & 1.857 & 2.611\\
 &  & Ref.~\cite{woowang1989} & 1.078 & 1.216 & 1.286 & 1.396 & 1.574 & 1.892 & 2.498\\
\cline{2-10} 
\multirow{6}{*}{$\dfrac{\pi}{6}$} &  & Polygonal SBFEM ($p=$ 2) & 0.738 & 0.871 & 0.947 & 1.037 & 1.153 & 1.315 & 1.554\tabularnewline
 & $F_{I}$ & Ref.~\cite{dauxmoes2000} & - & 0.864 & 0.943 & 1.031 & 1.145 & 1.297 & 1.531\\
 &  & Ref.~\cite{woowang1989} & 0.730 & 0.872 & 0.948 & 1.038 & 1.155 & 1.315 & 1.543\\
\cline{3-10} 
 &  & Polygonal SBFEM ($p=$ 2) & 0.153 & 0.342 & 0.419 & 0.464 & 0.504 & 0.559 & 0.637\\
 & $F_{II}$ & Ref.~\cite{dauxmoes2000} & - & 0.351 & 0.421 & 0.463 & 0.504 & 0.556 & 0.631\\
 &  & Ref.~\cite{woowang1989} & 0.141 & 0.336 & 0.417 & 0.463 & 0.504 & 0.556 & 0.639\\
\cline{2-10} 
\multirow{6}{*}{$\dfrac{\pi}{3}$} &  & Polygonal SBFEM ($p=$ 2) & 0.069 & 0.162 & 0.220 & 0.264 & 0.306 & 0.348 & 0.393\tabularnewline
 & $F_{I}$ & Ref.~\cite{dauxmoes2000} & - & 0.155 & 0.212 & 0.260 & 0.301 & 0.348 & 0.388\\
 &  & Ref.~\cite{woowang1989} & 0.039 & 0.147 & 0.210 & 0.257 & 0.300 & 0.344 & 0.390\\
\cline{3-10} 
 &  & Polygonal SBFEM ($p=$ 2) & 0.158 & 0.345 & 0.425 & 0.469 & 0.502 & 0.534 & 0.566\\
 & $F_{II}$ & Ref.~\cite{dauxmoes2000} & - & 0.347 & 0.424 & 0.466 & 0.500 & 0.532 & 0.566\\
 &  & Ref.~\cite{woowang1989} & 0.137 & 0.332 & 0.418 & 0.465 & 0.500 & 0.532 & 0.565\\
\hline 
\end{tabular}
\label{tab:twocrkholeresu}
\end{table}


\section{Conclusions}
In this paper, we compared three different displacement based finite element formulations over arbitrary polygons. The accuracy and the convergence properties of different approaches were studied with a few benchmark problems from linear elasticity. From the numerical examples, it is seen that the Polygonal SBFEM yields the most accurate results. Although, the Polygonal FEM yields optimal results, the numerical integration should be carried out accurately to improve the results, as discussed in \cite{talischipaulino2013}. The results from the Polygonal nSFEM are stiffer and yields sub-optimal convergence in both the displacement norm and the energy norm. It is noted that for this study, we had employed cell-based smoothed FEM, whilst other smoothing techniques are possible. The Polygonal SBFEM eliminates the need to compute the shape functions or to employ special numerical integration technique. Only the boundary of the polygon is discretized and the stiffness matrix is computed directly. Moreover, an extension to higher order shape functions along the boundary of the polygon is straightforward. When applied to linear elastic fracture mechanics, the scaled boundary polygon formulation yields accurate results. It is noted that the polygonal SBFEM does not require enrichment of the approximation basis and the stress intensity factors can be computed directly. This approach can readily be combined with the extended FEM. By combining the scaled boundary formulation with the XFEM, a priori knowledge of the asymptotic fields is not required.

\section*{Acknowledgements} 
S Natarajan would like to acknowledge the financial support of the School of Civil and Environmental Engineering, The University of New South Wales for his research fellowship for the period Sep. 2012 onwards. 


\bibliographystyle{spphys}       
\bibliography{psbfem}   

\end{document}